\documentclass[11pt,a4paper]{article}
\usepackage{amsfonts}
\usepackage{latexsym}
\usepackage{cite}
\usepackage{amsmath,amsfonts,latexsym,amssymb}
\usepackage[mathscr]{eucal}
\usepackage{cases,color}
\usepackage{amsthm}
\usepackage{hyperref}
\usepackage[bf,small]{caption2}
\usepackage{float}
\usepackage{graphicx}
\usepackage{amsmath}
\usepackage{amssymb}
\usepackage[all]{xy}

\newtheorem{theorem}{theorem}[section]

\newtheorem{cor}[theorem]{Corollary}
\newtheorem{conv}[theorem]{Convention}

\newtheorem{exmp}[theorem]{Example}
\newtheorem{lem}[theorem]{Lemma}

\newtheorem{que}[theorem]{Question}
\newtheorem{rmk}[theorem]{Remark}
\newtheorem{thm}[theorem]{Theorem}

\begin{document}

\title{\vspace{-1.5cm}\textbf{The ${\rm SL}(2,\mathbb{C})$-character variety of an arborescent knot}}
\author{\Large Haimiao Chen}

\date{}
\maketitle

\begin{abstract}
  We describe a procedure for computing the ${\rm SL}(2,\mathbb{C})$-character variety of an arborescent knot. Along the way, we clarify several facts about representations of arborescent tangles. Then we study a family of hyperbolic knots whose exteriors contain closed essential surfaces, showing that each of these knots has $1$-dimensional character variety.
  This provides infinitely many positive answers to a question of Boyer and Zhang posed in 1998.

\end{abstract}

\section{Introduction}

Fix $G={\rm SL}(2,\mathbb{C})$ throughout. Given a finitely presented group $\Gamma$, a homomorphism $\rho:\Gamma\to G$ is called a
{\it $G$-representation} of $\Gamma$. Call $\rho$ reducible if the elements of ${\rm Im}(\rho)$ have a common eigenvector, and irreducible otherwise. Call $\mathcal{R}(\Gamma):=\hom(\Gamma,G)$ the {\it $G$-representation variety} of $\Gamma$.
For each $\rho\in\mathcal{R}(\Gamma)$, its {\it character} is defined as the function $\chi_\rho:\Gamma\to\mathbb{C}$ sending $x$ to
${\rm tr}(\rho(x))$. It is known that two irreducible representations $\rho,\rho'$ have the same character if and only if they are conjugate, meaning that there exists $\mathbf{a}\in G$ such that $\rho'(x)=\mathbf{a}\rho(x)\mathbf{a}^{-1}$ for all $x\in\Gamma$.
The set $\mathcal{X}(\Gamma)=\{\chi_\rho\colon\rho\in\mathcal{R}(\Gamma)\}$ turns out to be an algebraic set, and is called the
$G$-{\it character variety} of $\Gamma$. The subset $\mathcal{X}^{\rm irr}(\Gamma)$ of $\mathcal{X}(\Gamma)$ consisting of characters of irreducible representations is Zariski open, and is called the {\it irreducible character variety}.

For a $3$-manifold $M$, abbreviate $\mathcal{X}(\pi_1(M))$ to $\mathcal{X}(M)$.
For a link $L\subset S^3$, let $\pi(L)=\pi_1(E_L)$, where $E_L$ denotes the exterior of $L$, abbreviate $\mathcal{R}(\pi(L))$ to $\mathcal{R}(L)$ and call it the $G$-representation variety of $L$, and so forth.

Nowadays, character varieties play a significant role in low-dimensional topology, but are still considered notoriously difficult to compute.
Very little is known about character varieties of general links. So it is worthwhile to describe $\mathcal{X}(L)$ for as many links $L$ as possible.
On the other hand, we believe it hopeful to reveal structural properties for certain families of knots.

In \cite{Ch22} the author proposed a method for determining $\mathcal{X}^{\rm irr}(K)$ for each Montesinos knot $K$, and showed a decomposition $\mathcal{X}^{\rm irr}(K)=\mathcal{X}_0\sqcup\mathcal{X}_1\sqcup\mathcal{X}_2\sqcup\mathcal{X}'$, such that each part has a distinguished feature.
In this paper, we aim to extend the scope to arborescent knots. These knots form an interesting class, and had been studied extensively by Bonahon and Siebenmann \cite{BS16}, and are abundant among those with small number of crossings. In particular, among the $84$ knots with at most $9$ crossings in Rolfsen's table, only $8_{18},9_{34},9_{39},9_{40},9_{41},9_{47},9_{49}$ are non-arborescent.
As another important fact, almost all arborescent knots are hyperbolic \cite{Wu11}.

It can be expected that new phenomena will be uncovered by studying knots other than ``small" knots such as torus knots, 2-bridge knots and 3-strand pretzel knots, as usually treated in the literature.

Here by studying a family of arborescent knots, we positively answer the following question of Boyer and Zhang, posed as \cite[Question 10.3]{BZ98}.
\begin{que}\label{que:BZ98}
Is there a hyperbolic knot whose exterior contains a closed essential surface but the ${\rm PSL}(2,\mathbb{C})$-character variety is 1-dimensional?
\end{que}
This is the first time to find infinitely many hyperbolic knots whose exteriors contain essential surfaces not detected by ideal points of character varieties. See Section \ref{sec:conclusion} for detailed discussions.

The content is organized as follows. In Section 2, we present some preliminaries on $2\times 2$ matrices and tangles, as well as representations of tangles. In Section 3, we describe a procedure of computing $\mathcal{X}^{\rm irr}(K)$ for arborescent knots $K$. Most of the work is concerned with representations of arborescent tangles, about which several important facts are clarified.
In Section 4, we study a family of hyperbolic arborescent knots, showing that $\dim\mathcal{X}(K)=1$ for each $K$ in that family. It is known that $E_K$ contains closed essential surfaces, thus $K$ fulfills the required conditions in Question \ref{que:BZ98}.

\section{Preparations}

\subsection{Algebraic aspects}

We always use bold letters to denote elements of $G$.

Let $\mathsf{T}_+,\mathsf{T}_-,\mathsf{T}_0$ respectively denote the subset of $G$ consisting of upper-triangular, lower-triangular, diagonal matrices.

Let $\mathbf{e}$ denote the $2\times 2$ identity matrix.
Let
\begin{alignat*}{2}
%&\mathbf{w}=\left(\begin{array}{cc} 0 & 1 \\ -1 & 0 \end{array}\right), \
\mathbf{p}&=\left(\begin{array}{cc} 1 & 1 \\ 0 & 1 \end{array}\right), \qquad
&\mathbf{d}(\kappa)&=\left(\begin{array}{cc} \kappa & 0 \\ 0 & \kappa^{-1} \end{array}\right), \\
\mathbf{u}_{\kappa}(\xi)&=\left(\begin{array}{cc} \kappa & \xi \\ 0 & \kappa^{-1} \end{array}\right), \qquad
&\mathbf{l}_{\kappa}(\xi)&=\left(\begin{array}{cc} \kappa & 0 \\ \xi & \kappa^{-1}\end{array}\right).
\end{alignat*}

For $t\in\mathbb{C}$, let $G(t)=\{\mathbf{x}\in G\setminus\{\pm\mathbf{e}\}\colon {\rm tr}(\mathbf{x})=t\}$.

For $\mathbf{a},\mathbf{b}\in G$, let $\mathbf{a}\lrcorner\mathbf{b}=\mathbf{a}\mathbf{b}\mathbf{a}^{-1}$.

Call a tuple $(\mathbf{a}_1,\ldots,\mathbf{a}_k)\in G^k$ {\it reducible} if $\mathbf{a}_1,\ldots,\mathbf{a}_k$ share an eigenvector, and
{\it irreducible} otherwise. Say that $(\mathbf{a}'_1,\ldots,\mathbf{a}'_k)$ is {\it conjugate} to $(\mathbf{a}_1,\ldots,\mathbf{a}_k)$ if there exists $\mathbf{c}\in G$ with $\mathbf{c}\lrcorner\mathbf{a}_i=\mathbf{a}'_i$ for $1\le i\le k$.
Call a reducible pair $(\mathbf{a}_1,\mathbf{a}_2)$ {\it nonabelian reducible} (NR for short) if $\mathbf{a}_1\mathbf{a}_2\ne\mathbf{a}_2\mathbf{a}_1$.

For $\lambda\ne -1$ and $\mu\ne 0$, put
\begin{align*}
\mathbf{h}_{t}^{\lambda}(\mu)=
\frac{1}{\lambda+1}\left(\begin{array}{cc} \lambda t & \mu \\ (t^2-\lambda-\lambda^{-1}-2)\lambda\mu^{-1} & t \end{array}\right).
\end{align*}
For $t\ne 0$, put
\begin{align*}
\mathbf{k}_{t}(\alpha)=\left(\begin{array}{cc} \alpha+t/2 & (t^2/4-1-\alpha^2)/(2t) \\ 2t & -\alpha+t/2 \end{array}\right).
\end{align*}
By direct computation,
\begin{align}
{\rm tr}\big(\mathbf{h}_{t}^{\lambda}(\mu)^{-1}\mathbf{h}_{t}^{\lambda}(\nu)\big)
=\frac{2t^2+(\lambda+\lambda^{-1}+2-t^2)(\mu\nu^{-1}+\mu^{-1}\nu)}{\lambda+\lambda^{-1}+2}.   \label{eq:tr-h}
\end{align}

For the following lemma, refer to \cite[Lemma 2.2]{Ch22}:
\begin{lem}\label{lem:key}
Let $\mathbf{a}_1,\mathbf{a}_2\in G(t)$.
\begin{enumerate}
  \item[\rm(a)] If $\mathbf{a}_1\mathbf{a}_2=\mathbf{d}(\lambda)$ with $\lambda+\lambda^{-1}\ne \pm2, t^2-2$,
                then $\mathbf{a}_1=\mathbf{h}_{t}^{\lambda}(-\lambda\mu)$, $\mathbf{a}_2=\mathbf{h}_{t}^{\lambda}(\mu)$ for some $\mu\ne 0$.
  \item[\rm(b)] If $\mathbf{a}_1\mathbf{a}_2=\mathbf{p}$, then $\mathbf{a}_1, \mathbf{a}_2\in\mathsf{T}_+$;
                more precisely, $\mathbf{a}_1=\mathbf{u}_{\kappa^{-1}}(\xi)$,
                $\mathbf{a}_2=\mathbf{u}_{\kappa}(\kappa-\xi)$ for some $\xi$ and $\kappa$ with $\kappa+\kappa^{-1}=t$.
  \item[\rm(c)] If $\mathbf{a}_1\mathbf{a}_2=-\mathbf{p}$ and $t\ne 0$, then $\mathbf{a}_1=\mathbf{k}_{t}(\alpha+t)$,
                $\mathbf{a}_2=\mathbf{k}_{t}(\alpha)$ for some $\alpha$.
  \item[\rm(d)] If $\mathbf{a}_1\mathbf{a}_2=\mathbf{d}(\kappa^2)$ with $\kappa+\kappa^{-1}=t\ne0,\pm 2$, then one of the following cases occurs:
                {\rm(i)} $\mathbf{a}_1=\mathbf{a}_2=\mathbf{d}(\kappa)$; {\rm(ii)} $\mathbf{a}_1=\mathbf{u}_{\kappa}(\kappa^2\xi)$, $\mathbf{a}_2=\mathbf{u}_{\kappa}(-\xi)$ for some $\xi\ne0$; {\rm(iii)} $\mathbf{a}_1=\mathbf{l}_{\kappa}(-\xi)$, $\mathbf{a}_2=\mathbf{l}_{\kappa}(\kappa^2\xi)$ for some $\xi\ne0$.
  \item[\rm(e)] The pair $(\mathbf{a}_1,\mathbf{a}_2)$ is reducible if and only if ${\rm tr}(\mathbf{a}_1\mathbf{a}_2)\in\{2,t^2-2\}$.
\end{enumerate}
\end{lem}

In the case $\mathbf{a}_1\mathbf{a}_2=\mathbf{p}$, the original statement in \cite{Ch22} is refined here.

Whenever ${\rm tr}(\mathbf{a}_1\mathbf{a}_2)=\lambda+\lambda^{-1}\ne\pm 2,t^2-2$, up to conjugacy we may assume $\mathbf{a}_1\mathbf{a}_2=\mathbf{d}(\lambda)$, so that $\mathbf{a}_1=\mathbf{h}_{t}^{\lambda}(-\lambda\mu)$ and $\mathbf{a}_2=\mathbf{h}_{t}^{\lambda}(\mu)$ for some $\mu\ne 0$. This will often simplify computations.

To illustrate the convenience, we establish a useful lemma:
\begin{lem}\label{lem:regular}
Suppose $\mathbf{a}_1,\mathbf{a}_2,\mathbf{a}_3,\mathbf{a}_4\in G(t)$ satisfy $\mathbf{a}_1\mathbf{a}_2=\mathbf{a}_3\mathbf{a}_4=\mathbf{c}$, with
${\rm tr}(\mathbf{c})\ne 2,t^2-2$.
\begin{enumerate}
  \item[\rm(i)] If ${\rm tr}(\mathbf{a}_1\mathbf{a}_3^{-1})=2$, then $\mathbf{a}_1=\mathbf{a}_3$ and $\mathbf{a}_2=\mathbf{a}_4$.
  \item[\rm(ii)] If $t\ne 0$, then $\mathbf{c}\ne-\mathbf{e}$, $\mathbf{a}_1\ne\mathbf{a}_3^{-1}$ and $\mathbf{a}_2\ne\mathbf{a}_4^{-1}$.
\end{enumerate}
\end{lem}

\begin{proof}
(i) Suppose ${\rm tr}(\mathbf{a}_1\mathbf{a}_3^{-1})=2$. Assume $\mathbf{a}_3^{-1}\mathbf{a}_1\ne\mathbf{e}$. Then up to conjugacy we may assume $\mathbf{a}_3^{-1}\mathbf{a}_1=\mathbf{p}$, so that also $\mathbf{a}_4\mathbf{a}_2^{-1}=\mathbf{p}$. By Lemma \ref{lem:key} (b), $\mathbf{a}_3,\mathbf{a}_1\in\mathsf{T}_+$, and $\mathbf{a}_4,\mathbf{a}_2\in\mathsf{T}_+$. This would imply
${\rm tr}(\mathbf{c})\in\{2,t^2-2\}$, contradicting the assumption. Thus, $\mathbf{a}_1=\mathbf{a}_3$, and $\mathbf{a}_2=\mathbf{a}_4$.

(ii) When ${\rm tr}(\mathbf{c})\ne -2$, up to conjugacy we may assume $\mathbf{c}=\mathbf{d}(\lambda)$ with $\lambda+\lambda^{-1}\ne\pm2,t^2-2$.
By Lemma \ref{lem:key} (a), $\mathbf{a}_2=\mathbf{h}^\lambda_t(\mu)$, $\mathbf{a}_4=\mathbf{h}^\lambda_t(\nu)$ for some $\mu,\nu$. Comparing the $(1,1)$-entries we see $\mathbf{a}_2\ne\mathbf{a}_4^{-1}$. So $\mathbf{a}_1\ne\mathbf{a}_3^{-1}$.

Now suppose ${\rm tr}(\mathbf{c})=-2$. If $\mathbf{c}=-\mathbf{e}$, then $\mathbf{a}_1=-\mathbf{a}_2^{-1}$, which would imply
$t={\rm tr}(\mathbf{a}_1)=-{\rm tr}(\mathbf{a}_2^{-1})=-t$, contradicting the assumption $t\ne 0$. Hence $\mathbf{c}\ne-\mathbf{e}$.
Up to conjugacy we may assume $\mathbf{c}=-\mathbf{p}$. By Lemma \ref{lem:key} (c), there exist $\alpha,\beta$ such that $\mathbf{a}_1=\mathbf{k}_{t}(\alpha)$, $\mathbf{a}_3=\mathbf{k}_{t}(\beta)$. Comparing the $(2,1)$-entries we immediately see $\mathbf{a}_1\ne\mathbf{a}_3^{-1}$, so that $\mathbf{a}_2\ne\mathbf{a}_4^{-1}$.
\end{proof}

Let
\begin{align*}
f_{t}(r_1,r_2,r_3;r)=r^2+t\left(t^2-\sum_{i=1}^3r_i\right)r+t^2\left(3-\sum_{i=1}^3r_i\right)+\sum_{i=1}^3r_i^2+r_1r_2r_3-4.   %\label{eq:f}
\end{align*}

\begin{lem}\label{lem:matrix}
Fix $t\in\mathbb{C}$.
\begin{enumerate}
  \item[\rm(a)] Given $t_{12}\ne 2,t^2-2$, up to conjugacy there exists a unique $(\mathbf{a}_1,\mathbf{a}_2)\in G(t)^2$
        with ${\rm tr}(\mathbf{a}_1\mathbf{a}_2)=t_{12}$.

        Given further $t_{13},t_{23},t_{123}$ with $f_{t}(t_{12},t_{13},t_{23};t_{123})=0$, there exists a unique $\mathbf{a}_3\in G(t)$ with
        $${\rm tr}(\mathbf{a}_1\mathbf{a}_3)=t_{13}, \qquad  {\rm tr}(\mathbf{a}_2\mathbf{a}_3)=t_{23}, \qquad
        {\rm tr}(\mathbf{a}_1\mathbf{a}_2\mathbf{a}_3)=t_{123}.$$
  \item[\rm(b)] Given any $t_{12},t_{13},t_{23},t_{123}$ with $f_{t}(t_{12},t_{13},t_{23};t_{123})=0$,
        there exists $(\mathbf{a}_1,\mathbf{a}_2,\mathbf{a}_3)\in G(t)^3$ such that
        $${\rm tr}(\mathbf{a}_{i_1}\cdots\mathbf{a}_{i_r})=t_{i_1\cdots i_r}, \qquad  1\le i_1<\cdots<i_r\le 3;$$
        when $(\mathbf{a}_1,\mathbf{a}_2,\mathbf{a}_3)$ is irreducible, it is unique up to conjugacy.
\end{enumerate}
\end{lem}
For a proof, see \cite{Go09}; in particular, Section 2 and Section 5.

\subsection{Combinatorial aspects}

By a tangle $T$ we simultaneously mean a tangle diagram and an embedded $1$-submanifold of a closed $3$-ball $B$ in $\mathbb{R}^3$ such that
$\partial B\cap T=\partial T$.

Given a tangle $T$, let $\mathcal{D}(T)$ denote the set of directed arcs.
For $\mathsf{a}\in\mathcal{D}(T)$, let $\mathsf{a}^{-1}$ denote the direct arc obtained from $\mathsf{a}$ by reversing the direction.
A map $\rho:\mathcal{D}(T)\to G$ is called a {\it representation} of $T$ if $\rho(\mathsf{a}^{-1})=\rho(\mathsf{a})^{-1}$ for all
$\mathsf{a}$ and $\rho(\mathsf{c})=\rho(\mathsf{a})\lrcorner\rho(\mathsf{b})$ for any $\mathsf{a},\mathsf{b},\mathsf{c}$ forming a crossing as
\begin{figure}[h]
  \centering
  \includegraphics[width=2.3cm]{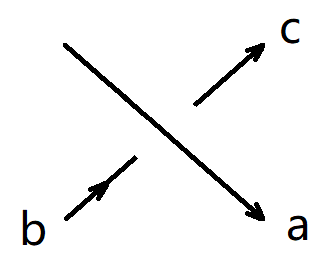}  \\
  %\caption{}\label{fig:crossing}
\end{figure}
\\
To present such $\rho$, it suffices to choose a direction for each arc and label the arc with an element of $G$.
In view of Wirtinger presentation, $\rho$ can be identified with a genuine representation $\tilde{\rho}:\pi_1(B\setminus T)\to G$.

Call $\rho$ reducible, irreducible, abelian, etc, if respectively $\tilde{\rho}$ is.
In particular, call $\rho$ NR if it is nonabelian and reducible.

Let $\mathcal{R}(T)$ denote the set of all representations of $T$.
Given $\rho\in\mathcal{R}(T)$ and $\mathbf{c}\in G$, let $\mathbf{c}\lrcorner\rho\in\mathcal{R}(T)$ send $\mathsf{a}$ to $\mathbf{c}\lrcorner\rho(\mathsf{a})$.
Call $\rho,\rho'\in\mathcal{R}(T)$ {\it conjugate} and denote $\rho\sim\rho'$ if $\rho'=\mathbf{c}\lrcorner\rho$ for some $\mathbf{c}$.
Let $\mathcal{X}(T)=\mathcal{R}(T)/\sim$, the set of conjugacy classes of representations of $T$; let $\mathcal{X}^{\rm red}(T)$,
$\mathcal{X}^{\rm irr}(T)$ respectively denote the set of conjugacy classes of reducible and irreducible representations.

If $S$ is a subtangle of $T$, and $\rho\in\mathcal{R}(T)$, then the composite map $\mathcal{D}(S)\to\mathcal{D}(T)\stackrel{\rho}\to G$ is a representation of $S$, which is denoted by $\rho|_{S}$.

\begin{figure}[h]
  \centering
  \includegraphics[width=12cm]{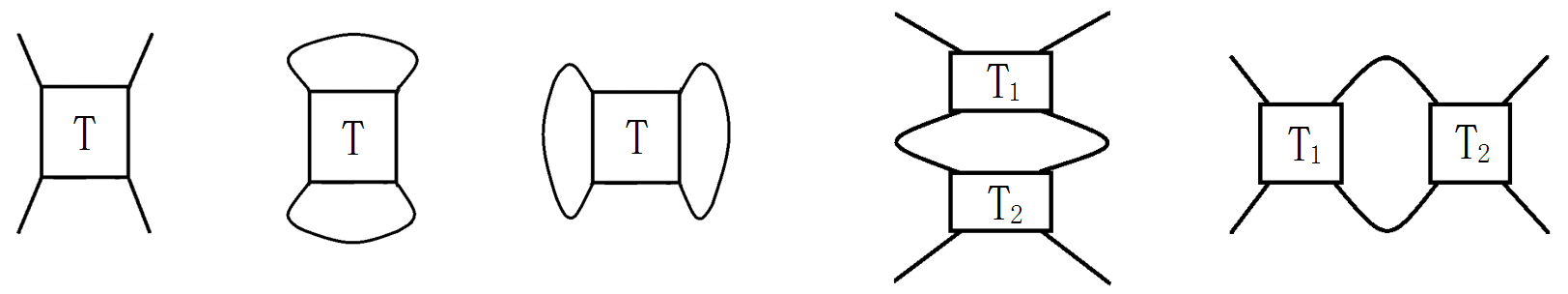}\\
  \caption{From left to right: a tangle $T$ in $\mathcal{T}_2^2$; $N(T)$; $D(T)$; $T_1\ast T_2$; $T_1+T_2$.}\label{fig:tangle}
\end{figure}

Let $\mathcal{T}_2^2$ denote the set of four-end tangles of the form shown at leftmost in Figure \ref{fig:tangle}.
To each $T\in\mathcal{T}_2^2$ are associated two links, called the {\it numerator} $N(T)$ and the {\it denominator} $D(T)$.
Defined on $\mathcal{T}_2^2$ are the vertical composition $\ast$ and the horizontal composition $+$.

\begin{figure}[h]
  \centering
  \includegraphics[width=4cm]{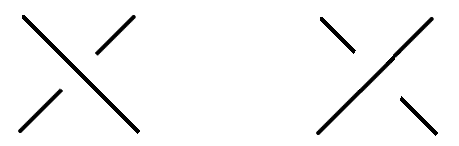}\\
  \caption{Left: the tangle $[1]$. Right: $[-1]$.}\label{fig:basic}
\end{figure}

A tangle $T\in\mathcal{T}_2^2$ is called {\it arborescent} if it can be constructed from copies of $[\pm1]$ (see Figure \ref{fig:basic}) via repeatedly applying $\ast$ and $+$. Let $\mathcal{T}_{\rm ar}$ denote the set of arborescent tangles.
An {\it arborescent knot} is a knot of the form $N(T)$ or $D(T)$ for some $T\in\mathcal{T}_{\rm ar}$.

For $k\ne 0$, the horizontal composite of $|k|$ copies of $[1]$ (resp. $[-1]$) is denoted by $[k]$ if $k>0$ (resp. $k<0$), and the
vertical composite of $|k|$ copies of $[1]$ (resp. $[-1]$) is denoted by $[1/k]$ if $k>0$ (resp. $k<0$).
Given $p/q\in\mathbb{Q}$, if a continued fraction of $p/q$ is
$$k_s+1/(k_{s-1}+\cdots+1/k_1)\cdots),$$
then the associated {\it rational tangle} is defined as
$$[p/q]=[[k_1],\ldots,[k_m]]=\begin{cases}
(\cdots([k_{1}]\ast[1/k_{2}])+\cdots)+[k_{m}], &2\nmid m \\
(\cdots([1/k_{1}]+[k_{2}])\ast\cdots)+[k_{m}], &2\mid m
\end{cases}.$$

When speaking of $T\in\mathcal{T}_{\rm ar}$, we usually assume that an order of repeated horizontal/vertical compositions of rational tangles has been chosen. Call the arborescent subtangles appearing in the compositions {\it subsequent to} $T$.
For example, if $R_1,\ldots,R_6$ are rational tangles, then
$$T=((R_1+R_2)+(R_3\ast R_4))\ast (R_5+R_6)\in\mathcal{T}_{\rm ar},$$
and $R_3$, $R_3\ast R_4$, $(R_1+R_2)+(R_3\ast R_4)$, $R_5$, etc. are subsequent to $T$.

\section{Methodology}

Suppose $T\in\mathcal{T}_{\rm ar}$.
Let $T^{{\rm nw}}$, $T^{{\rm ne}}$, $T^{{\rm sw}}$, $T^{{\rm se}}$ respectively denote the arcs at the northwest, northeast, southwest, southeast end of $T$, all directed outward.

For $t\in\mathbb{C}$, let
$$\mathcal{R}_t(T)=\{\rho\in\mathcal{R}(T)\colon\rho(\mathsf{a})\in G(t)\ \text{for\ all\ }\mathsf{a}\in\mathcal{D}(T)\}.$$
Given $\rho\in\mathcal{R}_t(T)$ and $\o\in\{{\rm nw},{\rm ne},{\rm sw},{\rm se}\}$, let $\mathbf{x}^{\o}=\mathbf{x}^{\o}_\rho=\rho(T^{\o})$.
Put
\begin{align*}
\mathbf{z}_\rho=\mathbf{x}^{{\rm nw}}\mathbf{x}^{{\rm ne}}, \qquad
\dot{\mathbf{z}}_\rho=\mathbf{x}^{{\rm ne}}\mathbf{x}^{{\rm se}}, \qquad
\grave{\mathbf{z}}_\rho=\mathbf{x}^{{\rm nw}}\mathbf{x}^{{\rm se}},  \qquad \acute{\mathbf{z}}_\rho&=\mathbf{x}^{{\rm sw}}\mathbf{x}^{{\rm ne}}.
\end{align*}
Furthermore, we associate
\begin{align*}
u&=u(\rho)={\rm tr}(\mathbf{z}_\rho), \qquad   \dot{u}=\dot{u}(\rho)={\rm tr}(\dot{\mathbf{z}}_\rho), \\
\grave{u}&=\grave{u}(\rho)={\rm tr}(\grave{\mathbf{z}}_\rho),  \qquad \acute{u}=\acute{u}(\rho)={\rm tr}(\acute{\mathbf{z}}_\rho),  \\
&\hspace{18mm}  \check{u}=\check{u}(\rho)=\grave{u}-\acute{u}.
\end{align*}

We mainly consider $\rho\in\mathcal{R}_t(T)$ for a fixed $t$.
The {\it trace-free} case, i.e. the one with $t=0$, had been thoroughly investigated in \cite{Ch19}.

\subsection{Representations of rational tangles revisited} \label{sec:rational}

It is worth clarifying some facts about representations of rational tangles.

\begin{figure}[h]
  \centering
  \includegraphics[width=11cm]{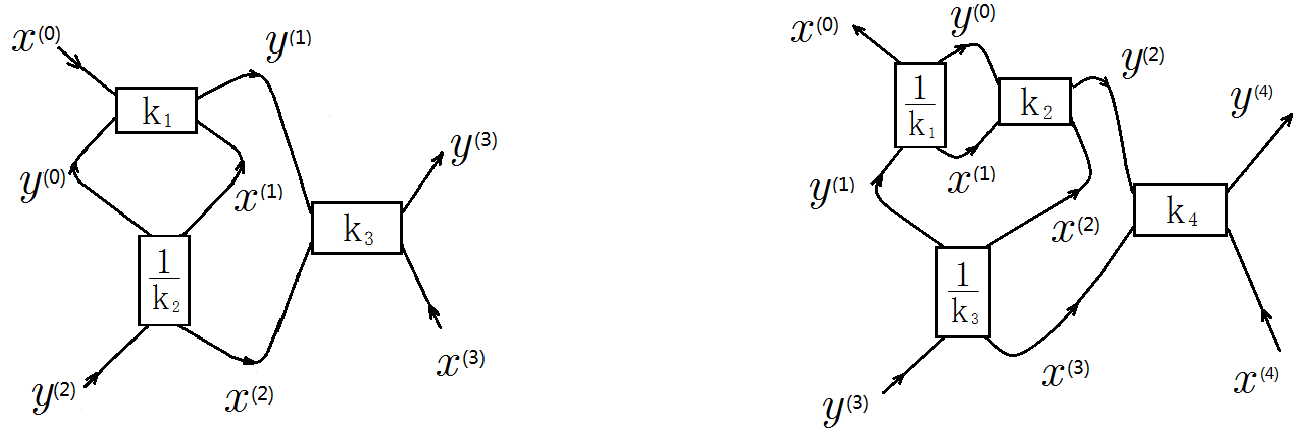}\\
  \caption{The rational tangle $[[k_1],\ldots,[k_m]]$; the case $m$ is odd (resp. even) is shown at left (resp. right).
  The $x^{(j)}$'s, $y^{(j)}$'s are directed arcs.}\label{fig:rational}
\end{figure}

In this subsection, suppose $T$ is a rational tangle, and let $x^{(j)},y^{(j)}\in\mathcal{D}(T)$ denote the directed arcs as shown in
Figure \ref{fig:rational}.

Given $\rho\in\mathcal{R}_t(T)$, let $\mathbf{x}^{(j)}=\rho(x^{(j)})$, $\mathbf{y}^{(j)}=\rho(y^{(j)})$;
let $\mathbf{x}=\mathbf{x}^{(0)}$, $\mathbf{y}=\mathbf{y}^{(0)}$, and $r={\rm tr}(\mathbf{x}\mathbf{y})$.
As explained in \cite[Section 3.1]{Ch22}, each of $\mathbf{x}^{(j)}$ and $\mathbf{y}^{(j)}$ can be written as a $\mathbb{Z}[t,r]$-linear combination of $\mathbf{e}$, $\mathbf{x}$, $\mathbf{y}$, $\mathbf{x}\mathbf{y}$.
Hence $\rho$ is determined by $\mathbf{x}$, $\mathbf{y}$; denote $\rho$ as $\rho^{\mathbf{x}}_{\mathbf{y}}$.
Clearly, $u(\rho),\dot{u}(\rho),\acute{u}(\rho),\grave{u}(\rho)$ are all determined by $t,r$. Let $u(t,r)=u_T(t,r)=u(\rho)\in\mathbb{Z}[t,r]$, etc.

\begin{rmk}\label{rmk:rep-rational}
\rm Recall \cite[Lemma 3.1]{Ch22} that $\mathbf{y}$ can be written as a word in
$\mathbf{x}^{\rm nw},\mathbf{x}^{\rm ne},\mathbf{x}^{\rm sw}, \mathbf{x}^{\rm se}$. Remembering $\mathbf{x}^{\rm nw}\in\{\mathbf{x}^{\pm1}\}$ and
$\mathbf{x}^{\rm se}=(\mathbf{x}^{\rm sw}\mathbf{x}^{\rm nw}\mathbf{x}^{\rm ne})^{-1}$, we see that $\rho$ is reducible if and only if $(\mathbf{x}^{\rm nw},\mathbf{x}^{\rm ne},\mathbf{x}^{\rm sw})$ is.
\end{rmk}

\begin{lem}\label{lem:rational-irreducible}
Suppose $r\ne 2,t^2-2$. Given $\mathbf{a}^{\rm nw}, \mathbf{a}^{\rm ne}, \mathbf{a}^{\rm sw}\in G(t)$ with
\begin{alignat*}{2}
{\rm tr}(\mathbf{a}^{\rm nw}\mathbf{a}^{\rm ne})&=u(t,r), \qquad  &&{\rm tr}(\mathbf{a}^{\rm nw}\mathbf{a}^{\rm sw})=\dot{u}(t,r),  \\
{\rm tr}(\mathbf{a}^{\rm sw}\mathbf{a}^{\rm ne})&=\acute{u}(t,r), \qquad  &&{\rm tr}(\mathbf{a}^{\rm sw}\mathbf{a}^{\rm nw}\mathbf{a}^{\rm ne})=t,
\end{alignat*}
there exists a unique pair $(\mathbf{x},\mathbf{y})\in G(t)^2$ such that ${\rm tr}(\mathbf{x}\mathbf{y})=r$ and
$\rho^{\mathbf{x}}_{\mathbf{y}}(T^{\o})=\mathbf{a}^{\o}$ for all $\o\in\{{\rm nw},{\rm ne},{\rm sw}\}$.
\end{lem}

\begin{proof}
By Lemma \ref{lem:matrix} (a), there exists an irreducible pair $(\mathbf{u},\mathbf{v})\in G(t)^2$ with
${\rm tr}(\mathbf{u}\mathbf{v})=r$. By Remark \ref{rmk:rep-rational},
$(\rho^{\mathbf{u}}_{\mathbf{v}}(T^{\rm nw}), \rho^{\mathbf{u}}_{\mathbf{v}}(T^{\rm ne}), \rho^{\mathbf{u}}_{\mathbf{v}}(T^{\rm sw}))$ is
irreducible. Due to the equalities on traces, there exists
$\mathbf{c}\in G$ such that
$\mathbf{c}\lrcorner\rho^{\mathbf{u}}_{\mathbf{v}}(T^{\o})=\mathbf{a}^{\o}$ for $\o\in\{{\rm nw},{\rm ne},{\rm sw}\}.$
Then $(\mathbf{x},\mathbf{y})=(\mathbf{c}\lrcorner\mathbf{u},\mathbf{c}\lrcorner\mathbf{v})$ fulfills the required conditions.

The uniqueness has been explained in Remark \ref{rmk:rep-rational}.
\end{proof}

\subsection{Representations of arborescent tangles}\label{sec:methodology}

The basic step is, for $T=T_1\ast T_2$ or $T=T_1+T_2$, to carefully analyze the relationship between $\mathcal{R}(T)$ and $\mathcal{R}(T_1),\mathcal{R}(T_2)$, in all possible cases.
The composition formulas in Lemma \ref{lem:composition} will be important.

By Lemma \ref{lem:regular} (ii), $\mathbf{z}=-\mathbf{e}$ is impossible, unless $t=0$.

When $t\ne 0$, in principle we can deal with the following cases separately: (i) $u\ne 2,t^2-2$; (ii) $\mathbf{z}=\mathbf{e}$;
(iii) $u=2$ but $\mathbf{z}\ne\mathbf{e}$, in which case, up to conjugacy we may assume $\mathbf{z}=\mathbf{p}$;
(iv) $u=t^2-2\ne \pm2$. In practice, each of the cases (ii)--(iv) has diverse possibilities and needs further investigation.
%When $\dot{u}\in\{2,t^2-2\}$ or $u(\rho)\in\{2,t^2-2\}$, we can apply Lemma \ref{lem:key} (b), (d).

We focus on (i), the generic case, and allow $t=0$.

When $u=\lambda+\lambda^{-1}\ne \pm 2,t^2-2$, up to conjugacy we may assume $\mathbf{z}=\mathbf{d}(\lambda)$, then by Lemma \ref{lem:key} (a),  there exist $\mu,\nu$ such that
\begin{alignat*}{2}
\mathbf{x}^{\rm nw}&=\mathbf{h}^\lambda_t(-\lambda\nu), \qquad  &\mathbf{x}^{\rm ne}&=\mathbf{h}^\lambda_t(\nu), \\
(\mathbf{x}^{\rm sw})^{-1}&=\mathbf{h}^\lambda_t(-\lambda\mu),  \qquad &(\mathbf{x}^{\rm se})^{-1}&=\mathbf{h}^\lambda_t(\mu).
\end{alignat*}
An application of (\ref{eq:tr-h}) leads to
\begin{align}
(u+2)\dot{u}&=2t^2+(u+2-t^2)(\nu\mu^{-1}+\nu^{-1}\mu), \label{eq:identity-0-1}  \\
(u+2)\grave{u}&=2t^2-(u+2-t^2)(\lambda\nu\mu^{-1}+\lambda^{-1}\nu^{-1}\mu),  \label{eq:identity-0-2}  \\
(u+2)\acute{u}&=2t^2-(u+2-t^2)(\lambda^{-1}\nu\mu^{-1}+\lambda\nu^{-1}\mu),  \label{eq:identity-0-3}  \\
(u+2)\check{u}&=(u+2-t^2)(\lambda-\lambda^{-1})(\nu^{-1}\mu-\nu\mu^{-1}).  \label{eq:identity-0-4}
\end{align}
The first three equations imply
\begin{align}
\grave{u}+\acute{u}+u\dot{u}=2t^2;  \label{eq:identity-1}
\end{align}
it follows from (\ref{eq:identity-0-1}), (\ref{eq:identity-0-4}) that
\begin{align}
\frac{\mu}{\nu}&=\frac{(u+2)(\dot{u}+\check{u}/(\lambda-\lambda^{-1}))-2t^2}{2(u+2-t^2)},  \label{eq:identity-2-1}  \\
\frac{\nu}{\mu}&=\frac{(u+2)(\dot{u}-\check{u}/(\lambda-\lambda^{-1}))-2t^2}{2(u+2-t^2)}.  \label{eq:identity-2-2}
\end{align}
The trivial identity $(\mu/\nu)\cdot(\nu/\mu)=1$ gives rise to
\begin{align}
\check{u}^2=(u-2)(\dot{u}-2)\big((u+2)(\dot{u}+2)-4t^2\big).  \label{eq:identity-3}
\end{align}
%which is equivalent to either of the following:
%\begin{align}
%u^2+\dot{u}^2+\grave{u}^2+u\dot{u}\grave{u}-2t^2(u+\dot{u}+\grave{u})+t^4+4t^2-4=0,  \label{eq:identity-4-1}  \\
%(\grave{u}-2)(\acute{u}-2)=(u+\dot{u}-t^2)^2.   \label{eq:identity-4-2}
%\end{align}

\begin{rmk}\label{rmk:identity}
\rm The equation (\ref{eq:identity-1}) implies that $u,\dot{u},\check{u}$ determine $\grave{u},\acute{u}$, so as to determine
$(\mathbf{x}^{\rm nw},\mathbf{x}^{\rm ne},\mathbf{x}^{\rm sw},\mathbf{x}^{\rm se})$ up to conjugacy.

Being polynomial identities, (\ref{eq:identity-1}), (\ref{eq:identity-3}) actually hold without the assumption $u\ne \pm 2,t^2-2$.
\end{rmk}

Define an equivalence relation $\sim$ on $\mathbb{C}^\ast\times\mathbb{C}^\ast$ by $(\mu,\nu)\sim(\mu^{-1},\nu^{-1})$.
Let
$$\mathcal{E}=\{(\mu,\nu)\in\mathbb{C}^\ast\times\mathbb{C}^\ast\colon\mu\ne\pm1\ \text{or\ }\nu\ne\pm1\}/\sim.$$
Partially defined on $\mathcal{E}$ are two operations:
\begin{alignat*}{2}
(\mu,\nu)\circ_h(\mu,\nu')&=(\mu,\nu\nu'), \qquad  &\mu\ne\pm1;  \\
(\mu,\nu)\circ_v(\mu',\nu)&=(\mu\mu',\nu), \qquad  &\nu\ne\pm1.
\end{alignat*}

When $u\ne \pm2,t^2-2$, define $\mathfrak{v}(\rho)=(\mu,\lambda)\in\mathcal{E}$ by taking $\lambda$ with $\lambda+\lambda^{-1}=u$ and putting
\begin{align}
\mu=\frac{(u+2)(\dot{u}+\check{u}/(\lambda-\lambda^{-1}))-2t^2}{2(u+2-t^2)}
\stackrel{(\ref{eq:identity-3})}=\Big(\frac{(u+2)(\dot{u}-\check{u}/(\lambda-\lambda^{-1}))-2t^2}{2(u+2-t^2)}\Big)^{-1}.  \label{eq:inherence-v}
\end{align}
Alternatively, up to conjugacy we may assume $\mathbf{z}=\mathbf{d}(\lambda)$, then by Lemma \ref{lem:key} (a),
\begin{alignat*}{2}
\mathbf{x}^{\rm nw}&=\mathbf{h}^\lambda_t(-\lambda\gamma), \qquad  &\mathbf{x}^{\rm ne}&=\mathbf{h}^\lambda_t(\gamma),  \\
(\mathbf{x}^{\rm sw})^{-1}&=\mathbf{h}^\lambda_t(-\lambda\gamma'), \qquad  &(\mathbf{x}^{\rm se})^{-1}&=\mathbf{h}^\lambda_t(\gamma')
\end{alignat*}
for some $\gamma,\gamma'$; we put $\mu=\gamma'/\gamma$.
It follows from (\ref{eq:inherence-v}) that
\begin{align}
\dot{u}&=\frac{2t^2+(u+2-t^2)(\mu+\mu^{-1})}{u+2}, \\
\check{u}&=\frac{(u+2-t^2)(\lambda-\lambda^{-1})(\mu-\mu^{-1})}{u+2}.
\end{align}

When $\dot{u}\ne \pm2,t^2-2$, define $\mathfrak{h}(\rho)=(\eta,\nu)\in\mathcal{E}$ by taking $\eta$ with $\eta+\eta^{-1}=\dot{u}$ and putting
\begin{align}
\nu=\frac{(\dot{u}+2)(u+\check{u}/(\eta-\eta^{-1}))-2t^2}{2(\dot{u}+2-t^2)}
\stackrel{(\ref{eq:identity-3})}=\Big(\frac{(\dot{u}+2)(u-\check{u}/(\eta-\eta^{-1}))-2t^2}{2(\dot{u}+2-t^2)}\Big)^{-1}.  \label{eq:inherence-h}
\end{align}
Alternatively, up to conjugacy we may assume $\dot{\mathbf{z}}=\mathbf{d}(\eta)$, then by Lemma \ref{lem:key} (a),
\begin{alignat*}{2}
\mathbf{x}^{\rm ne}&=\mathbf{h}^\eta_t(-\eta\gamma), \qquad   &\mathbf{x}^{\rm se}&=\mathbf{h}^\eta_t(\gamma),  \\
(\mathbf{x}^{\rm nw})^{-1}&=\mathbf{h}^\eta_t(-\eta\gamma'), \qquad &(\mathbf{x}^{\rm sw})^{-1}&=\mathbf{h}^\eta_t(\gamma')
\end{alignat*}
for some $\gamma,\gamma'$; we put $\nu=\gamma/\gamma'$.
It follows from (\ref{eq:inherence-h}) that
\begin{align}
u&=\frac{2t^2+(\dot{u}+2-t^2)(\nu+\nu^{-1})}{\dot{u}+2}, \\
\check{u}&=\frac{(\dot{u}+2-t^2)(\eta-\eta^{-1})(\nu-\nu^{-1})}{\dot{u}+2}.
\end{align}

By definition, the following properties are immediate:
\begin{align}
\mathfrak{h}(\rho_1+\rho_2)&=\mathfrak{h}(\rho_1)\circ_h\mathfrak{h}(\rho_2),  \label{eq:composition-h}  \\
\mathfrak{v}(\rho_1\ast\rho_2)&=\mathfrak{v}(\rho_1)\circ_v\mathfrak{v}(\rho_2).  \label{eq:composition-v}
\end{align}

\begin{exmp}\label{exmp:integer}
\rm Let $T$ be a rational tangle, and $\rho\in\mathcal{R}_t^{\rm irr}(T)$. Write $\dot{u}=\eta+\eta^{-1}$, $u=\mu+\mu^{-1}$.

When $T=[1]$, up to conjugacy we may assume $(\mathbf{x}^{\rm nw})^{-1}=\mathbf{h}_t^{\eta}(-\eta)$, $(\mathbf{x}^{\rm sw})^{-1}=\mathbf{h}_t^{\eta}(1)$, then $\mathbf{x}^{\rm se}=\mathbf{h}_t^{\eta}(-\eta)$, and it is easy to see $\mathfrak{h}(\rho)=(\eta,-\eta)$.
Alternatively, up to conjugacy we may also assume $(\mathbf{x}^{\rm sw})^{-1}=\mathbf{h}_t^{\mu}(-\mu)$,
$(\mathbf{x}^{\rm se})^{-1}=\mathbf{h}_t^{\mu}(1)$, then $\mathbf{x}^{\rm nw}=\mathbf{h}_t^{\mu}(1)$, and it is easy to see $\mathfrak{v}(\rho)=(-\mu,\mu)$.

For $r=\zeta+\zeta^{-1}$ and $n\in\mathbb{Z}$, put $\theta_n(r)=\zeta^n+\zeta^{-n}$,
\begin{align*}
\omega_n(r)&=\begin{cases} (\zeta^n-\zeta^{-n})/(\zeta-\zeta^{-1}), &\zeta\ne\pm 1 \\ n\zeta^{n-1}, &\zeta\in\{\pm 1\} \end{cases}, \\
\alpha_n(r)&=\frac{2t^2+(-1)^n(r+2-t^2)\theta_n(r)}{r+2}.
\end{align*}
In particular,
\begin{align*}
\alpha_{\pm2}(r)&=r^2-2+t^2(2-r)=2+(r-2)(r+2-t^2),  \\
\alpha_{\pm3}(r)&=3r-r^3+t^2(r-1)^2=2-(r+2-t^2)(r-1)^2.
\end{align*}
When $T=[k]$, repeatedly applying (\ref{eq:composition-h}), we obtain $\mathfrak{h}(\rho)=(\eta,(-\eta)^k)$, so
\begin{align}
u=\alpha_k(\dot{u}),  \qquad    \check{u}=(-1)^k(\dot{u}-2)(\dot{u}+2-t^2)\omega_{k}(\dot{u}).  \label{eq:[k]}
\end{align}
When $T=[1/k]$, repeatedly applying (\ref{eq:composition-v}), we obtain $\mathfrak{v}(\rho)=((-\mu)^k,\mu)$, so
\begin{align}
\dot{u}=\alpha_k(u),  \qquad    \check{u}=(-1)^k(u-2)(u+2-t^2)\omega_{k}(u).   \label{eq:[1/k]}
\end{align}

In general case, one may refer to the formulas given in \cite[Section 3.1]{Ch22}, which were expressed in terms of the polynomials $\psi_j(n)$.
\end{exmp}

\begin{rmk}\label{rmk:inherence-t=0}
\rm In the case $t=0$, things can be dramatically simplified. It is easy to see that $\mathfrak{h}(\rho)=\mathfrak{v}(\rho)=(\mu,\nu)$, with
$$\dot{u}=\mu+\mu^{-1},  \qquad  u=\nu+\nu^{-1},  \qquad  \check{u}=(\mu-\mu^{-1})(\nu-\nu^{-1}).$$

When $T=[p/q]=[[k_1],\ldots,[k_m]]$, as shown in \cite[Example 4.4]{Ch19}, up to conjugacy, each $\rho\in\mathcal{R}_0^{\rm irr}([p/q])$ is determined by $\mathfrak{v}(\rho)=\mathfrak{h}(\rho)=(\alpha,\beta)$ with
$$(-\alpha)^p+(-\beta)^q=0,$$
through $r=\alpha^{-\tilde{q}}\beta^{\tilde{p}}+\alpha^{\tilde{q}}\beta^{-\tilde{p}}$, where $\tilde{p}/\tilde{q}=[[k_2,\ldots,k_m]]^{(-1)^{m-1}}$.
\end{rmk}

\medskip

Given $b_1,c_1,b_2,c_2\in\mathbb{C}$ and $a\ne 2,t^2-2$, put
\begin{align}
g_a(b_1,c_1;b_2,c_2)&=\frac{(a+2)b_1b_2+2t^2(2-b_1-b_2)+c_1c_2/(a-2)}{2(a+2-t^2)},  \label{eq:g}  \\
h_a(b_1,c_1;b_2,c_2)&=\frac{(a+2)(b_1c_2+b_2c_1)-2t^2(c_1+c_2)}{2(a+2-t^2)}.    \label{eq:h}
\end{align}

\begin{lem}\label{lem:composition}
Suppose $T=T_1\ast T_2$ or $T=T_1+T_2$. Given $\rho\in\mathcal{R}_t(T)$, let $u=u(\rho)$, $\dot{u}=\dot{u}(\rho)$, and let $u_i=u(\rho|_{T_i})$, $\dot{u}_i=\dot{u}(\rho|_{T_i})$, $\check{u}_i=\check{u}(\rho|_{T_i})$, $i=1,2$.
\begin{enumerate}
  \item[\rm(i)] If $T=T_1\ast T_2$ and $u\ne 2,t^2-2$, then
        \begin{align*}
        \dot{u}(\rho)=g_u(\dot{u}_1,\check{u}_1;\dot{u}_2,\check{u}_2), \qquad
        \check{u}(\rho)=h_u(\dot{u}_1,\check{u}_1;\dot{u}_2,\check{u}_2).
        \end{align*}
  \item[\rm(ii)] If $T=T_1+T_2$ and $\dot{u}\ne 2,t^2-2$, then
        \begin{align*}
        u(\rho)=g_{\dot{u}}(u_1,\check{u}_1;u_2,\check{u}_2),  \qquad
        \check{u}(\rho)=h_{\dot{u}}(u_1,\check{u}_1;u_2,\check{u}_2).
        \end{align*}
\end{enumerate}
\end{lem}

\begin{proof}
We only prove (i); the proof for (ii) is similar.

Suppose $u\ne -2$. Suppose $\mathfrak{v}(\rho|_{T_i})=(\mu_i,\lambda)$, $i=1,2$, with $\lambda+\lambda^{-1}=u$.
Then by (\ref{eq:inherence-v}),
\begin{align*}
\mu_i^{\pm1}=\frac{1}{2(u+2-t^2)}\left((u+2)\left(\dot{u}_i\pm\frac{\check{u}_i}{\lambda-\lambda^{-1}}\right)-2t^2\right),
\end{align*}
and by (\ref{eq:composition-v}), $\mathfrak{v}(\rho)=(\mu_1\mu_2,\lambda)$.
By direct computation,
\begin{align*}
\dot{u}&=\frac{2t^2+(u+2-t^2)(\mu_1\mu_2+\mu_1^{-1}\mu_2^{-1})}{u+2}=g_u(\dot{u}_1,\check{u}_1;\dot{u}_2,\check{u}_2),  \\
\check{u}&=\frac{(u+2-t^2)(\lambda-\lambda^{-1})(\mu_1^{-1}\mu_2^{-1}-\mu_1\mu_2)}{u+2}=h_u(\dot{u}_1,\check{u}_1;\dot{u}_2,\check{u}_2).
\end{align*}

By continuity, the assertion also holds when $u=-2$ (in which case $t\ne 0$, due to $u\ne t^2-2$).
\end{proof}

For $a,b,c\in\mathbb{C}$, let $\mathcal{X}_{t;a,b}^c(T)$ denote the set of conjugacy classes of representations $\rho\in\mathcal{R}_t(T)$ with $u(\rho)=a$, $\dot{u}(\rho)=b$, $\check{u}(\rho)=c$.

\begin{cor}\label{cor:composite}
Suppose $a,b,c\in\mathbb{C}$.
\begin{enumerate}
  \item[\rm(i)] When $a\ne 2,t^2-2$, the map $\rho\mapsto(\rho|_{T_1},\rho|_{T_2})$ induces a bijection
        $$\mathcal{X}_{t;a,b}^c(T_1\ast T_2) \cong{\bigcup}_{(b_1,c_1;b_2,c_2)\in g_a^{-1}(b)\cap h_a^{-1}(c)}\mathcal{X}_{t;a,b_1}^{c_1}(T_1)
        \times\mathcal{X}_{t;a,b_2}^{c_2}(T_2).$$
  \item[\rm(ii)] When $b\ne 2,t^2-2$, the map $\rho\mapsto(\rho|_{T_1},\rho|_{T_2})$ induces a bijection
        $$\mathcal{X}_{t;a,b}^c(T_1+T_2)\cong
        {\bigcup}_{(a_1,c_1;a_2,c_2)\in g_b^{-1}(a)\cap h_b^{-1}(c)}\mathcal{X}_{t;a_1,b}^{c_1}(T_1)
        \times\mathcal{X}_{t;a_2,b}^{c_2}(T_2).$$
\end{enumerate}
\end{cor}

\begin{proof}
We only prove (i); the proof for (ii) is similar.

Given $\rho_i\in\mathcal{R}_t(T_i)$ with $u(\rho_i)=a$, $\dot{u}(\rho_i)=b_i$, $\check{u}(\rho_i)=c_i$, $i=1,2$,
due to
$${\rm tr}(\rho_1(T_1^{\rm sw})\rho_1(T_1^{\rm se}))={\rm tr}(\rho_2(T_2^{\rm nw})^{-1}\rho_2(T_2^{\rm ne})^{-1})=a\ne 2,t^2-2,$$
by Lemma \ref{lem:matrix} (i), there exists $\mathbf{c}\in G$ such that
$\mathbf{c}\lrcorner\rho_1(T_1^{\rm sw})=\rho_2(T_2^{\rm nw})^{-1}$ and
$\mathbf{c}\lrcorner\rho_1(T_1^{\rm se})=\rho_2(T_2^{\rm ne})^{-1}.$
We can glue $\mathbf{c}\lrcorner\rho_1$ and $\rho_2$ into $\rho\in\mathcal{R}_t(T)$. By Lemma \ref{lem:composition} (i), $\dot{u}(\rho)=g_a(b_1,c_1;b_2,c_2)$, $\check{u}(\rho)=h_a(b_1,c_1;b_2,c_2)$.

In this way, we construct the inverse to the map
$$\mathcal{X}_{t;a,b}^c(T_1\ast T_2) \to{\bigcup}_{(b_1,c_1;b_2,c_2)\in g_a^{-1}(b)\cap h_a^{-1}(c)}\mathcal{X}_{t;a,b_1}^{c_1}(T_1)
\times\mathcal{X}_{t;a,b_2}^{c_2}(T_2)$$
induced by $\rho\mapsto(\rho|_{T_1},\rho|_{T_2})$.
\end{proof}

\subsection{Reducible representations}

When $S=S_1+S_2$ or $S=S_1\ast S_2$, possibly $\rho|_{S_1}$ or $\rho|_{S_2}$ is reducible for some $\rho\in\mathcal{R}^{\rm irr}(S)$.
So it is worth taking $\mathcal{R}^{\rm red}(T)$ into account for $T\in\mathcal{T}_{\rm ar}$.

Suppose $\rho\in\mathcal{R}_t(T)$ is NR, so that $t\ne \pm2$. Fix $\kappa\ne\pm1$ with $\kappa+\kappa^{-1}=t$.
Up to conjugacy we may assume ${\rm Im}(\rho)\subset\mathsf{T}_+$.
Assume that for $\mathsf{a}\in\mathcal{D}(T)$,
$$\rho(\mathsf{a})=\mathbf{u}_{\kappa}(\xi_{\mathsf{a}})^{\epsilon_{\mathsf{a}}}
=\mathbf{u}_{\kappa^{\epsilon_{\mathsf{a}}}}(\epsilon_{\mathsf{a}}\xi_{\mathsf{a}}), \qquad  \xi_{\mathsf{a}}\in\mathbb{C}, \quad  \epsilon_{\mathsf{a}}\in\{\pm1\}.$$

\begin{figure}[H]
  \centering
  \includegraphics[width=7cm]{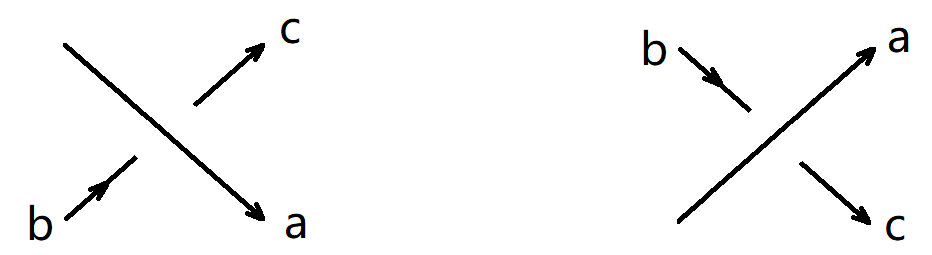}\\
  \caption{}\label{fig:pm}
\end{figure}

When $\mathsf{a}, \mathsf{b}, \mathsf{c}$ form a positive crossing (as shown in Figure \ref{fig:pm}, left), $\rho(\mathsf{c})=\rho(\mathsf{a})\lrcorner\rho(\mathsf{b})$, i.e.
$\mathbf{u}_{\kappa}(\xi_{\mathsf{c}})^{\epsilon_{\mathsf{c}}}
=\mathbf{u}_{\kappa}(\xi_{\mathsf{a}})^{\epsilon_{\mathsf{a}}}\lrcorner(\mathbf{u}_{\kappa}(\xi_{\mathsf{b}})^{\epsilon_{\mathsf{b}}})$.
Comparing the $(1,1)$- and $(1,2)$-entries, we obtain
\begin{align*}
\kappa^{\epsilon_{\mathsf{c}}}=\kappa^{\epsilon_{\mathsf{b}}},  \qquad
\epsilon_{\mathsf{c}}\xi_{\mathsf{c}}=\kappa^{2\epsilon_{\mathsf{a}}}\epsilon_{\mathsf{b}}\xi_{\mathsf{b}}
+\kappa^{\epsilon_{\mathsf{a}}}(\kappa^{-\epsilon_{\mathsf{b}}}-\kappa^{\epsilon_{\mathsf{b}}})\epsilon_{\mathsf{a}}\xi_{\mathsf{a}},
\end{align*}
which imply
\begin{align}
\xi_{\mathsf{c}}=\big(1-\kappa^{2\epsilon_{\mathsf{a}}}\big)\xi_{\mathsf{a}}+\kappa^{2\epsilon_{\mathsf{a}}}\xi_{\mathsf{b}}. \label{eq:relation-1}
\end{align}
Similarly, when $\mathsf{a}, \mathsf{b}, \mathsf{c}$ form a negative crossing, one has
\begin{align}
\xi_{\mathsf{c}}=\big(1-\kappa^{-2\epsilon_{\mathsf{a}}}\big)\xi_{\mathsf{a}}+\kappa^{-2\epsilon_{\mathsf{a}}}\xi_{\mathsf{b}}. \label{eq:relation-2}
\end{align}

Orient the components of $T$ (as a $1$-submanifold of a $3$-ball), so that a direction is chosen for each of its arcs. Numerate the crossings of $T$ as $\mathfrak{c}_1,\ldots,\mathfrak{c}_n$, and the directed arcs as $x_1,\ldots,x_{n+2}$.
Construct $Q_T=(q_{ij})_{n\times(n+2)}$ as follows: for each $i$, if $\mathfrak{c}_i$ is formed by $\mathsf{a}=x_{k_1}$, $\mathsf{b}=x_{k_2}$, $\mathsf{c}=x_{k_3}$, then let
\begin{align*}
q_{i,k_1}=1-\kappa^{2\varepsilon_i}, \qquad  q_{i,k_2}=\kappa^{2\varepsilon_i},   \qquad  q_{i,k_3}=-1,
\end{align*}
where $\varepsilon_i=\epsilon_{\mathsf{a}}$ (resp. $\varepsilon_i=-\epsilon_{\mathsf{a}}$) if $\mathfrak{c}_i$ is positive (resp. negative);
let $q_{i,j}=0$ for $j\notin\{k_1,k_2,k_3\}$.
By (\ref{eq:relation-1}), (\ref{eq:relation-2}), we have $Q_T(\xi_{x_1},\ldots,\xi_{x_{n+2}})=0.$

\begin{rmk}
\rm The matrix $Q_T$ identifies with the $Q_T$ defined in \cite[Section 4.2]{Ch26}, under which $\kappa^{2\varepsilon_i}$ corresponds to $\tau_{\overline{i}}^{\epsilon}$.
%In this way, the relationship between reducible ${\rm SL}(2,\mathbb{C})$-representations of $D(T)$ and the Alexander polynomial of $D(T)$ is clarified.
\end{rmk}

For $S\in\mathcal{T}_{\rm ar}$, define its {\it fraction} $f(S)\in\mathbb{Q}\cup\{\infty\}$ by setting $f(S)=p/q$ if $S=[p/q]$, and recursively,
$$f(S_1+S_2)=f(S_1)+f(S_2),  \qquad  \frac{1}{f(S_1\ast S_2)}=\frac{1}{f(S_1)}+\frac{1}{f(S_2)},$$
with the convention $\infty+\infty=\infty$, $1/0=\infty$, $1/\infty=0$.

Let $\xi^{\rm nw}=\xi_{T^{\rm nw}}$, etc.
As elucidated in \cite[Section 4]{Ch26}, if $T$ is generic, in the sense that $f(S)\ne 0,\infty$ for all tangles subsequent to $T$,
then for each arc $x$, there exists $b^x,c^x\in\mathbb{Q}(\kappa)$ such that
$$\xi_x=(1-b^x)\xi^{\rm nw}+b^x\xi^{\rm ne}=(1-c^x)\xi^{\rm nw}+c^x\xi^{\rm sw}.$$
In particular, there exist $b^{\rm sw},b^{\rm se},c^{\rm ne},c^{\rm se}\in\mathbb{Q}(\kappa)$ determined by $T$ such that
\begin{alignat}{2}
\left(\begin{array}{cc} \xi^{\rm sw} \\ \xi^{\rm se} \end{array}\right)&=F_v^T
\left(\begin{array}{cc} \xi^{\rm nw} \\ \xi^{\rm ne} \end{array}\right), \qquad
&&\text{with} \quad  F_v^T=\left(\begin{array}{cc} 1-b^{\rm sw} & b^{\rm sw} \\ 1-b^{\rm se} & b^{\rm se} \end{array}\right),  \label{eq:Fv} \\
\left(\begin{array}{cc} \xi^{\rm ne} \\ \xi^{\rm se} \end{array}\right)&=F_h^T
\left(\begin{array}{cc} \xi^{\rm nw} \\ \xi^{\rm sw} \end{array}\right), \qquad
&&\text{with} \quad  F_h^T=\left(\begin{array}{cc} 1-c^{\rm ne} & c^{\rm ne} \\ 1-c^{\rm se} & c^{\rm se} \end{array}\right).  \label{eq:Fh}
\end{alignat}
When $T$ is not generic, the situation is obscure, and more efforts are required.

When $\kappa^2=-1$, i.e. $t=0$, we can recursively show $b^{\rm sw}=-1/f(T)$, $c^{\rm ne}=-f(T)$.
The equations (\ref{eq:Fv}), (\ref{eq:Fh}) are consistent with the following result presented in \cite[Page 110 and 112]{Ch19}:
If $\rho\in\mathcal{R}_0^{\rm red}(T)$, so that $\dot{u}=2\alpha$, $u=2\beta $ with $\alpha,\beta\in\{\pm1\}$, then
\begin{align}
(\mathbf{x}^{\rm ne},\mathbf{x}^{\rm se})&=-\beta(\mathbf{x}^{\rm nw},\mathbf{x}^{\rm sw})C_{\alpha}(f(T)),   \label{eq:red-h}  \\
(\mathbf{x}^{\rm sw},\mathbf{x}^{\rm se})&=-\alpha(\mathbf{x}^{\rm nw},\mathbf{x}^{\rm ne})C_\beta(1/f(T)),   \label{eq:red-v}
\end{align}
where for $\gamma\in\{\pm1\}$ and $c\in\mathbb{C}$, we put
$$C_{\gamma}(c)=\left(\begin{array}{cc} 1+c  & -\gamma c  \\  \gamma c &  1-c \end{array}\right).$$

\subsection{On representations of arborescent knots}

%For $T\in\mathcal{T}_{\rm ar}$, t
The link $N(T)$ is equivalent to $D(T')$, where $T'$ is obtained from rotating $T$ by $\pi/2$.
Thus, we may just pay attention to knots of the form $K=D(T)$.
Evidently, representations of $K$ identify with $\rho\in\mathcal{R}(T)$ such that $\dot{\mathbf{z}}_\rho=\mathbf{e}$.

When $S=S_1+S_2$ (resp. $S=S_1\ast S_2$), call $\rho\in\mathcal{R}_t(S)$ non-degenerate if $\dot{u}(\rho)\ne 2,t^2-2$ (resp. $u(\rho)\ne 2,t^2-2$).
Call $\rho\in\mathcal{R}_t(T)$ {\it regular} if $\rho|_S$ is non-degenerate for each $S$ subsequent to $T$.
%that has the form $S=S_1+S_2$ or $S=S_1\ast S_2$.

Let $\mathcal{X}^{\rm reg}(K)\subset\mathcal{X}(K)$ denote the subset consisting of characters of regular representations.
Characters in $\mathcal{X}(K)\setminus\mathcal{X}^{\rm reg}(K)$ are relatively easy to determine, but the process may be tedious.

For $K=D(T)$ with $T=T_1\ast T_2$, we have a useful technique:
\begin{lem}\label{lem:convenient}
Suppose $T=T_1\ast T_2$, and $\rho\in\mathcal{R}_t(T)$ with $u(\rho)\ne 2,t^2-2$. Then $\dot{\mathbf{z}}_\rho=\mathbf{e}$ is equivalent to $\dot{u}(\rho|_{T_1})=\dot{u}(\rho|_{T_2})$ and $\check{u}(\rho|_{T_1})+\check{u}(\rho|_{T_2})=0$.
\end{lem}

\begin{proof}
If $\dot{\mathbf{z}}_\rho=\mathbf{e}$, then clearly,
$\dot{u}(\rho|_{T_1})=\dot{u}(\rho|_{T_2})$, $\acute{u}(\rho|_{T_1})=\grave{u}(\rho|_{T_2})$, and $\grave{u}(\rho|_{T_1})=\acute{u}(\rho|_{T_2});$
the latter two imply $\check{u}(\rho|_{T_1})+\check{u}(\rho|_{T_2})=0$.

Conversely, suppose $\dot{u}(\rho|_{T_1})=\dot{u}(\rho|_{T_2})$ and $\check{u}(\rho|_{T_1})+\check{u}(\rho|_{T_2})=0$. Let $a=u(\rho)=u(\rho|_1)=u(\rho|_2)$, $b=\dot{u}(\rho|_{T_1})$, $c=\check{u}(\rho|_{T_1})$.
By Lemma \ref{lem:composition} (i),
\begin{align*}
\dot{u}(\rho)=\ &g_a(b,c;b,-c)=\frac{(a+2)b^2+2t^2(2-2b)-c^2/(a-2)}{2(u+2-t^2)}  \\
\stackrel{(\ref{eq:identity-3})}=&\frac{(a+2)b^2+2t^2(2-2b)-(b-2)((a+2)(b+2)-4t^2)}{2(u+2-t^2)}=2.
\end{align*}
By Lemma \ref{lem:regular} (i), $\dot{\mathbf{z}}_\rho=\mathbf{e}$.
\end{proof}

\begin{rmk}
\rm As the proof shows, $\dot{\mathbf{z}}_\rho=\mathbf{e}$ is equivalent to the single condition $\dot{u}(\rho)=2$.
So actually there is redundance in $\dot{u}(\rho|_{T_1})=\dot{u}(\rho|_{T_2})$ and $\check{u}(\rho|_{T_1})+\check{u}(\rho|_{T_2})=0$. But the expression of $\dot{u}(\rho)$ is more complicated.
\end{rmk}

%{\color{red}Clearly write down the algorithm of computing $\mathcal{X}(T)$ step by step, as explicit as possible.
%Intensely aim to quickly write down the equations routinely.}

\begin{exmp}
\rm
\begin{figure}[h]
  \centering
  \includegraphics[width=11cm]{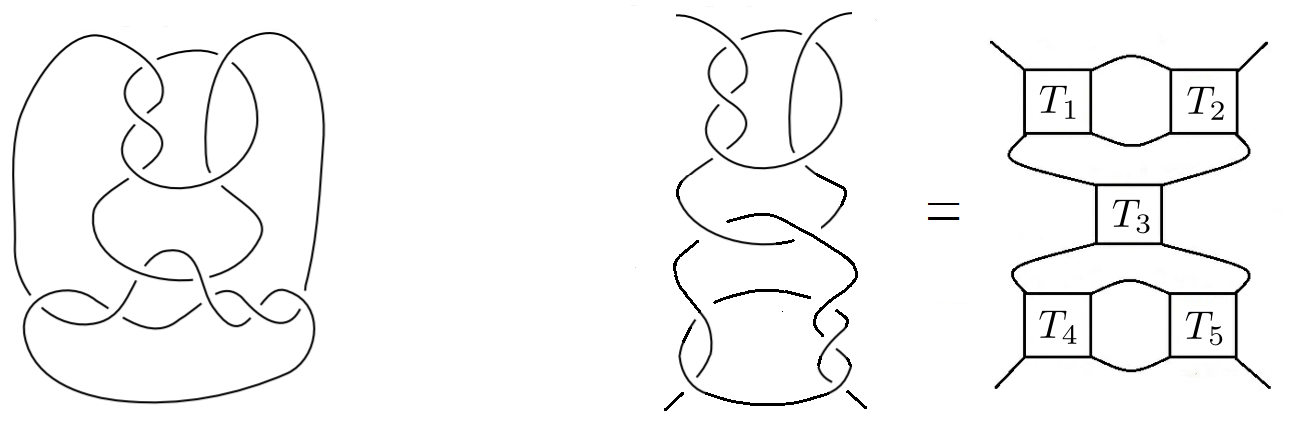}\\
  \caption{Left: the Conway knot $C$ (the figure is taken from \cite{Ho22}). Right: the tangle $T=(T_1+T_2)\ast T_3\ast(T_4+T_5)$, with $T_1=[1/3]$, $T_2=[-1/2]$, $T_3=[2]$, $T_4=[1/2]$, $T_5=[-1/3]$.}\label{fig:Conway}
\end{figure}

Let $C$ denote the famous Conway knot, shown in Figure \ref{fig:Conway}. It can be presented as $D(T)$ for
$T=(T_1+T_2)\ast T_3\ast(T_4+T_5),$
with
$T_1=[1/3]$, $T_2=[-1/2]$, $T_3=[2]$, $T_4=[1/2]$, $T_5=[-1/3].$
%Let $g\in\pi(C)$ denote the element represented by the dotted loop.

We shall determine $\mathcal{X}^{\rm reg}(C)$. The formulas (\ref{eq:[k]}), (\ref{eq:[1/k]}) will be used repeatedly, without more explanation.

Consider a general regular $\rho\in\mathcal{R}_t(T)$. Let $u_i=u(\rho|_{T_i})$, $\check{u}_{12}=\check{u}(\rho|_{T_1+T_2})$, etc.
Let $z={\rm tr}(\mathbf{z}_\rho)\ne 2,t^2-2$.

For the tangle $T_{12}:=T_1+T_2=[1/3]+[-1/2]$, let $x=\dot{u}_1=\dot{u}_2$, i.e.
$$x=3u_1-u_1^3+t^2(u_1-1)^2=u_2^2-2+t^2(2-u_2).$$
Note that
\begin{align*}
\check{u}_1&=(2-u_1)(u_1+2-t^2)(u_1^2-1)=2u_1^2+(x-2t^2)u_1+2t^2-4,  \\
\check{u}_2&=(2-u_2)(u_2+2-t^2)u_2=(2-x)u_2.
\end{align*}
When $x\ne 2,t^2-2$, applying Lemma \ref{lem:composition} to $T_1+T_2$, we can compute
\begin{align*}
z=u_{12}&=\frac{-u_1^2u_2+(t^2+1)u_1u_2-t^2u_1+(2-2t^2)u_2+2t^2}{x+2-t^2},     \\
\check{u}_{12}&=2u_1u_2-(x+2)z+2t^2(1-w),
\end{align*}
with
\begin{align*}
w=\frac{u_1^2-u_1u_2+(x+1-t^2)u_1+2u_2-2}{x+2-t^2}.
\end{align*}

Consider $T_{123}:=T_{12}\ast T_3=T_{12}\ast[2]$. Let $v=\dot{u}([2])$.
We have
\begin{align*}
z=u([2])&=v^2-2+t^2(2-v),  \\
\check{u}([2])&=(v-2)(v+2-t^2)v=(z-2)v.
\end{align*}
Applying Lemma \ref{lem:composition} to $T_{12}\ast T_3$, we can compute
\begin{align*}
\dot{u}_{123}&=\frac{v(u_1u_2+x-z-t^2w)+t^2(2-x)}{z+2-t^2},   \\
\check{u}_{123}&=\frac{(z+2)v(u_1u_2-x-z-t^2w)-t^2(2u_1u_2-(x+2)z+2t^2(1-w)-4v)}{z+2-t^2}.
\end{align*}

For the tangle $T_{45}:=T_4+T_5=[1/2]+[-1/3]$, let $y=\dot{u}_4=\dot{u}_5$, i.e.
$$y=u_4^2-2+t^2(2-u_4)=3u_5-u_5^3+t^2(u_5-1)^2.$$
Note that
\begin{align*}
\check{u}_4&=(u_4-2)(u_4+2-t^2)u_4=(y-2)u_4,  \\
\check{u}_5&=(u_5-2)(u_5+2-t^2)(u_5^2-1)=(2t^2-y)u_5-2u_5^2+4-2t^2.
\end{align*}
When $y\ne 2,t^2-2$, applying Lemma \ref{lem:composition} to $T_4+T_5$, we can compute
\begin{align*}
z=u_{45}&=\frac{-u_4u_5^2+(t^2+1)u_4u_5+(2-2t^2)u_4-t^2u_5+2t^2}{y+2-t^2},   \\
\check{u}_{45}&=(y+2)z-2u_4u_5-2t^2+\frac{2t^2(u_5^2-u_4u_5+2u_4+(y+1-t^2)u_5-2)}{y+2-t^2}.
\end{align*}

Applying Lemma \ref{lem:convenient} to $T_{123}\ast T_{45}$, we obtain the equations defining $\mathcal{X}^{\rm reg}(C)$
(with $x,y,z\ne 2,t^2-2$):
\begin{align*}
3u_1-u_1^3+t^2(u_1-1)^2=x,  \\
u_2^2-2+t^2(2-u_2)=x,   \\
-u_1^2u_2+(t^2+1)u_1u_2-t^2u_1+(2-2t^2)u_2+2t^2=(x+2-t^2)z, \\
u_4^2-2+t^2(2-u_4)=y,  \\
3u_5-u_5^3+t^2(u_5-1)^2=y,  \\
-u_4u_5^2+(t^2+1)u_4u_5+(2-2t^2)u_4-t^2u_5+2t^2=(y+2-t^2)z,  \\
v^2-2+t^2(2-v)=z,    \\
u_1^2-u_1u_2+(x+1-t^2)u_1+2u_2-2=(x+2-t^2)w,  \\
v(u_1u_2+x-z-t^2w)+t^2(2-x)=(z+2-t^2)y,  \\
\frac{(z+2)v(u_1u_2-x-z-t^2w)-t^2(2u_1u_2-(x+2)z+2t^2(1-w)-4v)}{z+2-t^2}\ \ \ \ \  \\
+(y+2)z-2u_4u_5-2t^2+\frac{2t^2(u_5^2-u_4u_5+2u_4+(y+1-t^2)u_5-2)}{y+2-t^2}=0.
\end{align*}
\end{exmp}

\section{A family of non-Montesinos knots}

\begin{figure}[H]
  \centering
  \includegraphics[width=12.8cm]{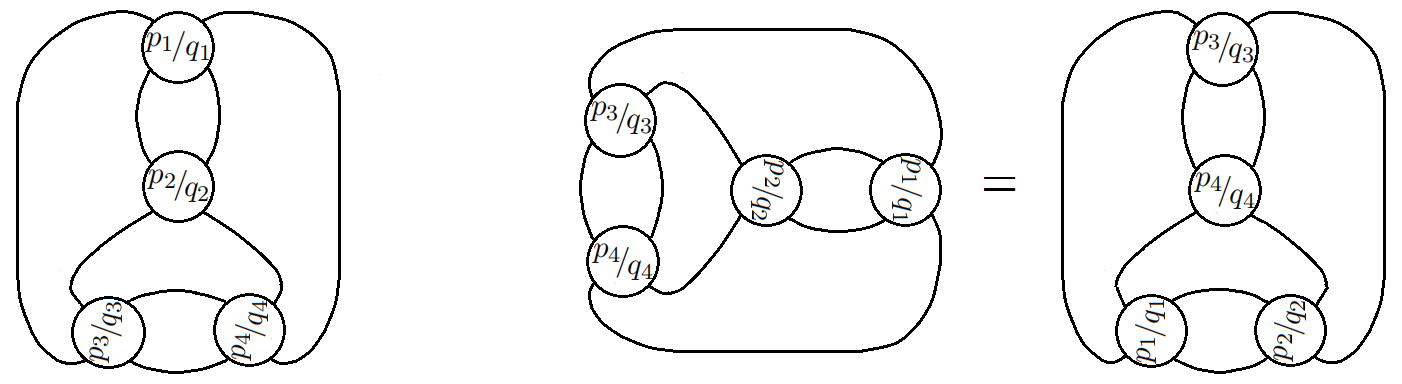}\\
  \caption{A non-Montesinos $3$-bridge knot.}\label{fig:3-bridge}
\end{figure}

By \cite[Theorem 1.1]{Ja11}, each non-Montesinos 3-bridge link falls into 3 families, one of which has the form
$K=D(T)$, with $T=T_1\ast T_2\ast(T_3+T_4)$, where
$$T_1=[p_1/q_1], \qquad  T_2=[p_2/q_2], \qquad  T_3=[-q_3/p_3], \qquad  T_4=[-q_4/p_4],$$
as depicted in the left part of Figure \ref{fig:3-bridge}.

From the right part of Figure \ref{fig:3-bridge} we see that $K$ is equivalent to $D(T')$, with $T'=T'_3\ast T'_4\ast(T'_1+T'_2)$, where $$T'_1=[-q_1/p_1], \qquad  T'_2=[-q_2/p_2], \qquad  T'_3=[p_3/q_3], \qquad  T'_4=[p_4/q_4].$$
It will be helpful to bear in the mind such a ``structural symmetry".

To ensure $K$ to be a knot, we assume
\begin{align}
2\nmid p_1p_2p_3p_4\ \ \ \text{or\ \ \ }2\nmid (p_1-p_2)(p_3-p_4).    \label{eq:assumption-0}
\end{align}
We further assume the following generic conditions:
\begin{align}
p_1,p_2,p_3,p_4\ge 3, \qquad (p_1,p_2)=(p_3,p_4)=1;    \label{eq:assumption-1}   \\
\frac{p_1}{q_1}+\frac{p_2}{q_2}\ne 0,    \quad    \frac{p_3}{q_3}+\frac{p_4}{q_4}\ne 0,  \qquad   \left(\frac{q_1}{p_1}+\frac{q_2}{p_2}\right)\left(\frac{q_3}{p_3}+\frac{q_4}{p_4}\right)\ne -1.   \label{eq:assumption-2}
\end{align}

The main goal of this section is to show
\begin{thm}\label{thm:dim=1}
Under the assumptions {\rm(\ref{eq:assumption-0})}--{\rm(\ref{eq:assumption-2})}, $\dim\mathcal{X}(K)=1$.
\end{thm}

We shall outline the steps of computing $\mathcal{X}^{\rm irr}(K)$, but do not present full details, to keep the paper reasonably short.

Consider a general element of $\mathcal{R}_t^{\rm irr}(K)$, which is the same as $\rho\in\mathcal{R}_t^{\rm irr}(T)$ with $\dot{\mathbf{z}}_\rho=\mathbf{e}$.
Let $\mathbf{z}=\mathbf{z}_\rho$, $z={\rm tr}(\mathbf{z})$, and let $\rho_1=\rho|_{T_1}$, $\rho_{34}=\rho|_{T_3+T_4}$, $\dot{\mathbf{z}}_{12}=\dot{\mathbf{z}}_{\rho_{12}}$, etc. Always remember $\mathbf{w}:=\dot{\mathbf{z}}_{34}=(\dot{\mathbf{z}}_{12})^{-1}$, so that $\dot{u}_{12}=\dot{u}_{34}=:w$.

Trace-free characters can be handled by the method developed in \cite{Ch19}.

When $t\ne 0$, up to conjugacy, we may assume one of the following holds:
$z\ne 2,t^2-2$; $\mathbf{z}=\mathbf{e}$; $\mathbf{z}=\mathbf{p}$;
$\mathbf{z}=\mathbf{d}(\kappa^2)$ with $\kappa+\kappa^{-1}=t\ne\pm2$.
Accordingly, $\mathcal{X}^{\rm irr}(K)\setminus\mathcal{X}_0^{\rm irr}(K)$ is decomposed into $4$ parts.
%For each part, the focus will be put on its dimension.

Due to $p_1/q_1+p_2/q_2\ne 0$, by \cite[Proposition 3.1]{CO15}, $L_{12}:=D(T_1\ast T_2)$ is a rational link (rather than an unlink). Similarly, $p_3/q_3+p_4/q_4\ne 0$ ensures $L_{34}:=N(T_3+T_4)$ to be rational link.
By \cite[Section 7]{ORS08}, if $L$ is a $n$-component rational link, $n\in\{1,2\}$, then $\dim\mathcal{X}(L)=n$;
actually, $\mathcal{X}(L)$ is locally parameterized by the value(s) taken at the meridian(s).
Thus, for each $t$, both of $\mathcal{X}_t(L_{12})$, $\mathcal{X}_t(L_{34})$ are finite sets.

\begin{conv}
\rm Recall the notations introduced in Section \ref{sec:rational}. Abbreviate $u_{T_i}(t,r_i)$, $\dot{u}_{T_i}(t,r_i)$, $\check{u}_{T_i}(t,r_i)$ respectively to $u_i(r_i)$, $\dot{u}_i(r_i)$, $\check{u}_i(r_i)$, omitting $t$, although they are polynomials in $t$ and $r_i$.

By Lemma \ref{lem:rational-irreducible}, to determine $\tau\in\mathcal{R}^{\rm irr}(T_i)$, it is sufficient to give
$\tau(T_i^{\rm nw})$, $\tau(T_i^{\rm ne})$, $\tau(T_i^{\rm sw})$ satisfying certain conditions.

We shall use the structural symmetry of $K$ to simplify the discussion.

In Section 4.2--4.5, we assume $t\ne 0$. Then by Lemma \ref{lem:regular} (ii), $\mathbf{z},\mathbf{w}\ne-\mathbf{e}$.
\end{conv}

\subsection{Trace-free representations}\label{sec:t=0}

Recall Remark \ref{rmk:inherence-t=0} that up to conjugacy each $\tau\in\mathcal{R}_0^{\rm irr}([p/q])$ is determined by $\mathfrak{v}(\tau)=\mathfrak{h}(\tau)=(\alpha,\beta)$ with
$(-\alpha)^p+(-\beta)^q=0$.

Throughout this subsection, suppose $\rho\in\mathcal{R}_0^{\rm irr}(K)$.

\begin{lem}\label{lem:t=0}
There are exactly three cases: {\rm(i)} $z,w\ne\pm 2$; {\rm(ii)} $z\ne\pm2$ and $\mathbf{w}\in\{\pm\mathbf{e}\}$;
{\rm(iii)} $\mathbf{z}\in\{\pm\mathbf{e}\}$ and $w\ne\pm2$.
%$\mathbf{g},\mathbf{w}\ne\pm\mathbf{p}$, and $\{\mathbf{g},\mathbf{w}\}\not\subseteq\{\pm\mathbf{e}\}$.
\end{lem}

\begin{proof}
It suffices to show that $\mathbf{z},\mathbf{w}\notin\{\pm\mathbf{p}\}$ and $\{\mathbf{z},\mathbf{w}\}\not\subseteq\{\pm\mathbf{e}\}$.

Assume $\mathbf{z}=\varepsilon\mathbf{p}$ with $\varepsilon\in\{\pm1\}$.
Then for $i=1,2$, one has
$\mathbf{x}_i^{\rm nw}(\mathbf{x}_i^{\rm ne})^{\varepsilon}=\mathbf{p}$, so $\mathbf{x}_i^{\rm nw},\mathbf{x}_i^{\rm ne}\in\mathsf{T}_+$ and
$\mathbf{x}_i^{\rm nw}\mathbf{x}_i^{\rm ne}\ne \mathbf{x}_i^{\rm ne}\mathbf{x}_i^{\rm nw}$;
by Lemma \ref{lem:key} (b), ${\rm Im}(\rho_i)\subset\mathsf{T}_+$.
Consequently, $w={\rm tr}(\mathbf{x}_1^{\rm nw}\mathbf{x}_2^{\rm sw})\in\{\pm2\}$. There are two possibilities:
\begin{enumerate}
  \item If $\mathbf{w}\in\{\pm\mathbf{e}\}$, then there exist $\alpha,\beta,\gamma\in\{\pm1\}$ such that
        $$\alpha(\mathbf{x}_1^{\rm nw},\mathbf{x}_1^{\rm ne})=(\mathbf{x}_2^{\rm sw},\mathbf{x}_2^{\rm se})
        =-\beta(\mathbf{x}_1^{\rm nw},\mathbf{x}_1^{\rm ne})C_\gamma(q_1/p_1+q_2/p_2);$$
        since $q_1/p_1+q_2/p_2\ne 0$, this implies that $\mathbf{x}_1^{\rm nw},\mathbf{x}_1^{\rm ne}$ are linearly dependent, contradicting $\mathbf{x}_1^{\rm nw}\mathbf{x}_1^{\rm ne}=\mathbf{z}=\varepsilon\mathbf{p}$.
  \item If $(\mathbf{x}_3^{\rm nw},\mathbf{x}_3^{\rm sw})$ is NR, then $(\nu\mathbf{d}(\eta))\lrcorner\mathbf{w}=\mathbf{p}$ for
        some $\nu\in\{\pm1\}$ and $\eta\ne 0$. Hence ${\rm Im}(\rho_3), {\rm Im}(\rho_4)\subset\mathsf{T}_+$.
        This contradicts the irreducibility of $\rho$.
\end{enumerate}

By symmetry, $\mathbf{w}\in\{\pm\mathbf{p}\}$ is neither possible.

Finally, we rule out $\{\mathbf{z},\mathbf{w}\}\subseteq\{\pm\mathbf{e}\}$. Assume on the contrary that $\mathbf{z}=\alpha\mathbf{e}$, $\mathbf{w}=\beta\mathbf{e}$ with $\alpha,\beta\in\{\pm1\}$. Then $\dot{u}_1=\beta\dot{u}_2$.
If $\rho_1$ is NR, then by (\ref{eq:red-h}),
$$(\mathbf{x}_1^{\rm ne},\mathbf{x}_1^{\rm se})=-\alpha(\mathbf{x}_1^{\rm nw},\mathbf{x}_1^{\rm sw})C_\gamma(p_1)
=(\mathbf{x}_1^{\rm ne},\mathbf{x}_1^{\rm se})C_\gamma(p_1),$$
which is impossible. Similarly, $\rho_2$ is not NR.
It follows that $\mathfrak{v}(\rho_1)=(\gamma,\alpha)$ and
$\mathfrak{v}(\rho_2)=(\beta/\gamma,\alpha)$ for some $\gamma$ with $\gamma+\gamma^{-1}=\dot{u}_1$.
We have
$$(-\gamma)^{p_1}+(-\alpha)^{q_1}=(-\beta/\gamma)^{p_2}+(-\alpha)^{q_2}=0.$$
Since $(p_1,p_2)=1$, we deduce $\gamma\in\{\pm1\}$. So $\rho_1$, $\rho_2$ are abelian.
By symmetry, $\rho_3,\rho_4$ are also abelian. But this contradicts the irreducibility of $\rho$.

\end{proof}

To investigate Case (i) in detail, suppose $z=\lambda+\lambda^{-1}$ and $w=\eta+\eta^{-1}$, with $\lambda,\eta\ne\pm1$.
For $i=1,2$, let $\mathfrak{v}(\rho_i)=(\alpha_i,\lambda)$, so that $\alpha_i+\alpha_i^{-1}=\dot{u}_i$; for $i=3,4$, let $\mathfrak{h}(\rho_i)=(\eta,\beta_i)$, so that $\beta_i+\beta_i^{-1}=u_i$.
Then $\lambda,\eta,\alpha_1,\alpha_2,\beta_3,\beta_4$ determine a representation if and only if
\begin{alignat}{3}
\alpha_1\alpha_2&=\eta,  \qquad  &(-\alpha_1)^{p_1}+(-\lambda)^{q_1}&=0,  \qquad  &(-\alpha_2)^{p_2}+(-\lambda)^{q_2}&=0,  \label{eq:tracefree-1}  \\
\beta_3\beta_4&=\lambda, \qquad  &(-\eta)^{-q_3}+(-\beta_3)^{p_3}&=0,  \qquad  &(-\eta)^{-q_4}+(-\beta_4)^{p_4}&=0.     \label{eq:tracefree-2}
\end{alignat}
Write $-\lambda=re^{i\phi}$ with $r>0$ and $0\le\phi<2\pi$. The last two equations in (\ref{eq:tracefree-1}) become
$$\alpha_1=-r^{\frac{q_1}{p_1}}e^{i\theta_1}, \qquad  \alpha_2=-r^{\frac{q_2}{p_2}}e^{i\theta_2}, \qquad \text{with\ \ \ }
\theta_j=\frac{q_j\phi+(2k_j+1)\pi}{p_j},$$
for some $k_1,k_2\in\mathbb{Z}$ with $0\le k_i<p_i$.
Consequently,
\begin{align}
\eta=\alpha_1\alpha_2=r^{\frac{q_1}{p_1}+\frac{q_2}{p_2}}e^{i(\theta_1+\theta_2)}.   \label{eq:eta}
\end{align}
Then (\ref{eq:tracefree-2}) can be converted into
\begin{align*}
r&=r^{-(\frac{q_3}{p_3}+\frac{q_4}{p_4})(\frac{q_1}{p_1}+\frac{q_2}{p_2})},  \\
\phi&\equiv -\left(\frac{q_3}{p_3}+\frac{q_4}{p_4}\right)(\theta_1+\theta_2+\pi)+\left(\frac{2k_3+1}{p_3}+\frac{2k_4+1}{p_4}-1\right)\pi\pmod{2\pi},
\end{align*}
for some $k_3,k_4\in\mathbb{Z}$ with $0\le k_i<p_i$. From (\ref{eq:assumption-2}) we see that $r=1$, and $\phi$ has finitely many choices.

For Case (ii), require $\eta\in\{\pm1\}$ and replace (\ref{eq:tracefree-2}) by
\begin{align}
(-\eta)^{-q_3}+(-\beta_3)^{p_3}=0,  \qquad  (-\eta)^{-q_4}+(-\beta_4)^{p_4}=0.     \label{eq:tracefree-2'}
\end{align}
Then $\beta_3,\beta_4$ has finitely many solutions. By (\ref{eq:eta}), $r=1$, and there exists $k\in\mathbb{Z}$ such that
$$\phi=\left(\frac{q_1}{p_1}+\frac{q_2}{p_2}\right)^{-1}\left(2k+\frac{1-\eta}{2}-\frac{2k_1+1}{p_1}-\frac{2k_2+1}{p_2}\right)\pi,$$
which determines $\alpha_1,\alpha_2$.
Up to conjugacy, we may fix $\mathbf{x}_i^{\rm nw}$ for $i=1,2,3$ and $\mathbf{x}_i^{\rm ne}$ for $i=1,2,4$.
Furthermore, by Lemma \ref{lem:matrix} (a), $\mathbf{x}_3^{\rm ne}$ is determined by $u_3,u_4$ and
$r:={\rm tr}(\mathbf{x}_3^{\rm ne}\mathbf{x}_3^{\rm nw}\mathbf{x}_4^{\rm ne})$, subject to $f_0(u_3,u_4,z;r)=0$.

\subsection{$z\ne 2,t^2-2$}\label{sec:generic}

Recall the functions introduced in (\ref{eq:g}), (\ref{eq:h}).

The condition $z\ne 2,t^2-2$ ensures $\rho_1,\rho_2$ to be irreducible.
We have
\begin{align}
u_1(r_1)=z, \qquad  u_2(r_2)=z,  \qquad  \dot{u}_3(r_3)=w, \qquad \dot{u}_4(r_4)=w,  \label{eq:Case1-1}  \\
g_z(\dot{u}_1(r_1),\check{u}_1(r_1);\dot{u}_2(r_2),\check{u}_2(r_2))=w.  \label{eq:Case1-2}
\end{align}
Given $t,z$, up to conjugacy, we may fix $\mathbf{x}_1^{\rm nw}=(\mathbf{x}_3^{\rm sw})^{-1}$,
$\mathbf{x}_1^{\rm ne}$, $\mathbf{x}_1^{\rm sw}=(\mathbf{x}_2^{\rm nw})^{-1}$,
$\mathbf{x}_1^{\rm se}=(\mathbf{x}_2^{\rm ne})^{-1}$, and then determine $\mathbf{x}_2^{\rm sw}=(\mathbf{x}_3^{\rm nw})^{-1}\in G(t)$ via
$${\rm tr}(\mathbf{x}_2^{\rm sw}\mathbf{x}_2^{\rm nw})=\dot{u}_2(r_2), \qquad
{\rm tr}(\mathbf{x}_2^{\rm sw}\mathbf{x}_2^{\rm ne})=\acute{u}_2(r_2), \qquad
{\rm tr}(\mathbf{x}_2^{\rm sw}\mathbf{x}_2^{\rm nw}\mathbf{x}_2^{\rm ne})=t,$$
so as to determine $\mathbf{x}_2^{\rm se}=(\mathbf{x}_4^{\rm ne})^{-1}$.
%By Lemma \ref{lem:rational-irreducible}, $\mathbf{x}_2^{\rm nw}$, $\mathbf{x}_2^{\rm ne}$, $\mathbf{x}_2^{\rm sw}$ indeed give rise to an irreducible representation of $T_2$.

\begin{enumerate}
  \item If $w\ne 2,t^2-2$, then $\rho_3,\rho_4$ are irreducible, and
        \begin{align}
        h_z(\dot{u}_1(r_1),\check{u}_1(r_1);\dot{u}_2(r_2),\check{u}_2(r_2))+h_w(u_3(r_3),\check{u}_3(r_3);u_4(r_4),\check{u}_4(r_4))=0.  \label{eq:Case1-3}
        \end{align}
        %, subject to (\ref{eq:Case1-1}), (\ref{eq:Case1-2}), (\ref{eq:Case1-3}).
        We can convert (\ref{eq:Case1-2}), (\ref{eq:Case1-3}) into polynomial equations, which together with (\ref{eq:Case1-1}) form a system of polynomial equations in $r_1,\ldots,r_4,w$; let $R(t,z)$ denote the resultant. A necessary condition for (\ref{eq:Case1-1})--(\ref{eq:Case1-3}) to have a solution is $R(t,z)=0$.

        Although we assume $t\ne 0$, the equations (\ref{eq:Case1-1})--(\ref{eq:Case1-3}) are also valid when $t=0$, in which case they are equivalent to (\ref{eq:tracefree-1}), (\ref{eq:tracefree-2}). From the result of Case (i) in Section \ref{sec:t=0} we see $R(0,z)\not\equiv0$.
        Hence
        $$\mathcal{V}:=\{(t,z)\in\mathbb{C}^2\colon R(t,z)=0\}$$
        is an algebraic curve. For each $(t,z)\in\mathcal{V}$, there exist finitely many pairs $(r_1,r_2)$ with $u_1(r_1)=u_2(r_2)=z$, and then $w$ is given by (\ref{eq:Case1-2}). Finally, there are finitely many pairs $(r_3,r_4)$ with $\dot{u}_3(r_3)=\dot{u}_4(r_4)=w$.
        Now $\rho$ is determined by $\mathbf{x}_3^{\rm ne}$, which is in turn determined by
        $${\rm tr}(\mathbf{x}_3^{\rm ne}\mathbf{x}_3^{\rm nw})=u_3(r_3), \qquad  {\rm tr}(\mathbf{x}_3^{\rm ne}\mathbf{x}_3^{\rm sw})=\acute{u}_3(r_3), \qquad  {\rm tr}(\mathbf{x}_3^{\rm ne}\mathbf{x}_3^{\rm sw}\mathbf{x}_3^{\rm nw})=t.$$
  \item When $w\in\{2,t^2-2\}$, let $R'(t,z)$ denote the resultant of the equations (\ref{eq:Case1-1}), (\ref{eq:Case1-2}),
        which also hold for $t=0$.
        From the result of Case (ii) in Section \ref{sec:t=0} we see $R'(t,z)\not\equiv 0$, so
        $$\mathcal{V}':=\{(t,z)\in\mathbb{C}^2\colon R'(t,z)=0\}$$
        is an algebraic curve. For each $(t,z)\in\mathcal{V}'$, up to finite ambiguity $r_1,r_2,r_3,r_4$ are determined via (\ref{eq:Case1-1}).
        Then $\rho$ is determined by $\mathbf{x}_3^{\rm ne}$ via ${\rm tr}(\mathbf{x}_3^{\rm ne}\mathbf{x}_3^{\rm nw})=u_3(r_3)$,
        ${\rm tr}(\mathbf{x}_3^{\rm ne}\mathbf{x}_4^{\rm ne})=t^2-u_4(r_4)$ and
        $r_5:={\rm tr}(\mathbf{x}_3^{\rm nw}\mathbf{x}_4^{\rm ne}\mathbf{x}_3^{\rm ne})$, subject to
        $$f_t\big(u_3(r_3),t^2-u_4(r_4),z;r_5\big)=0.$$

        Observe that when $w=2$, by Lemma \ref{lem:regular} (i), actually $\mathbf{w}=\mathbf{e}$.
\end{enumerate}

Therefore, the part of $\mathcal{X}^{\rm irr}(K)$ with $z\ne 2,t^2-2$ has dimension $1$.

\subsection{$\mathbf{z}=\mathbf{e}$}\label{sec:g=e}

In this case, $\mathbf{x}_1^{\rm ne}=(\mathbf{x}_1^{\rm nw})^{-1}$, $\mathbf{x}_2^{\rm ne}=(\mathbf{x}_2^{\rm nw})^{-1}$,
$\mathbf{x}_4^{\rm ne}=(\mathbf{x}_3^{\rm nw})^{-1}$.

For each $t$, there are finitely many pairs $(r_1,r_2)$ with $u_1(r_1)=u_2(r_2)=2$.

If $w\ne 2,t^2-2$, then the situation is symmetric to the subcase of Case 2 with $w=2$ in Section \ref{sec:generic}.

Suppose $w\in\{2,t^2-2\}$. Let $\epsilon=1$ if $w=2$, and $\epsilon=-1$ if $w=t^2-2$.
\begin{enumerate}
  \item Suppose $\mathbf{x}_1^{\rm nw}\mathbf{x}_3^{\rm nw}=\mathbf{x}_3^{\rm nw}\mathbf{x}_1^{\rm nw}$, i.e.
        $\mathbf{x}_3^{\rm nw}=(\mathbf{x}_1^{\rm nw})^{\epsilon}$.

        If $\rho_1,\rho_2$ are both abelian,
        then $\rho$ is determined by $\rho_{34}$, which identifies with a representation of $L_{34}$.

        Assume $\rho_1$ or $\rho_2$ is nonabelian. Consider the equations
        \begin{align}
        u_1(r_1)=u_2(r_2)=2, \qquad  \dot{u}_1(r_1)=w+\epsilon(\dot{u}_2(r_2)-2),  \label{eq:z=e}
        \end{align}
        which are also valid when $t=0$.
        Let $R^\ast(t)$ denote the resultant of (\ref{eq:z=e}) with the constraint $r_1,r_2\notin\{2,t^2-2\}$. 
        Precisely, let $x_1,\ldots,x_m$ be the roots of $u_1(r_1)-2$ other than $2,t^2-2$ in the algebraic closure of $\mathbb{Q}(t)$, and let $y_1,\ldots,y_n$ be those of $u_2(r_2)-2$; put
        $$R^\ast(t)=\prod_{i=1}^m\prod_{j=1}^n\big(\dot{u}_1(x_i)-w-\epsilon(\dot{u}_2(y_j)-2)\big).$$
        The proof of Lemma \ref{lem:t=0} shows that $z,w\in\{\pm2\}$ is impossible if $t=0$. Thus $R^\ast(0)\ne 0$. Only when $R^\ast(t)=0$, there exist $r_1,r_2\notin\{2,t^2-2\}$ satisfying (\ref{eq:z=e}), so that $\rho_1,\rho_2$ are irreducible.
        Then $\rho$ is determined by $(\mathbf{x}_1^{\rm nw},\mathbf{x}_2^{\rm nw},\mathbf{x}_3^{\rm ne})$, which up to conjugacy is determined by
        \begin{alignat*}{2}
        {\rm tr}(\mathbf{x}_1^{\rm nw}\mathbf{x}_2^{\rm nw})&=t^2-\dot{u}_1(r_1), \qquad
        &{\rm tr}(\mathbf{x}_1^{\rm nw}\mathbf{x}_3^{\rm ne})&=w+\epsilon(u_3(r_3)-2), \\
        {\rm tr}(\mathbf{x}_2^{\rm nw}\mathbf{x}_3^{\rm ne})&=r, \qquad
        &{\rm tr}(\mathbf{x}_1^{\rm nw}\mathbf{x}_2^{\rm nw}\mathbf{x}_3^{\rm ne})&=s,
        \end{alignat*}
        subject to $f_t\big(t^2-\dot{u}_1(r_1),w+\epsilon(u_3(r_3)-2),r;s\big)=0$.
  \item If $\mathbf{x}_1^{\rm nw}\mathbf{x}_3^{\rm nw}\ne\mathbf{x}_3^{\rm nw}\mathbf{x}_1^{\rm nw}$, then
        up to conjugacy we may assume $\mathbf{x}_1^{\rm nw}=\mathbf{d}(\kappa^\epsilon)$, $\mathbf{x}_3^{\rm nw}=\mathbf{u}_\kappa(1)$.
        For $\mathbf{x}_3^{\rm ne}$, there are two possibilities:
        \begin{itemize}
          \item If $\rho_3$ is irreducible, then $\mathbf{x}_3^{\rm ne}$ is determined by
                $${\rm tr}(\mathbf{x}_3^{\rm ne}\mathbf{x}_3^{\rm nw})=u_3(r_3), \quad
                {\rm tr}(\mathbf{x}_3^{\rm ne}\mathbf{x}_3^{\rm sw})=\acute{u}_3(r_3), \quad
                {\rm tr}(\mathbf{x}_3^{\rm ne}\mathbf{x}_3^{\rm sw}\mathbf{x}_3^{\rm nw})=t.$$
          \item If $\rho_3$ is NR, then $\mathbf{x}_3^{\rm ne}$ is determined by $\mathbf{x}_3^{\rm nw}$ and 
                $\mathbf{x}_3^{\rm sw}=(\mathbf{x}_1^{\rm nw})^{-1}$ through $F_h^{T_3}$ as in (\ref{eq:Fh}).
        %subject to $f_t(u_3(r_3),\acute{u}_3(r_3),\dot{u}_3(r_3),t)=0$.
        \end{itemize}
        For $\mathbf{x}_2^{\rm nw}$, there are three possibilities:
        \begin{itemize}
          \item If $\rho_1$ or $\rho_2$ is irreducible, then $\mathbf{x}_2^{\rm nw}$ is determined by $\dot{u}_1(r_1)$, $\dot{u}_2(r_2)$, $w$
                and $v:={\rm tr}(\mathbf{x}_1^{\rm nw}\mathbf{x}_2^{\rm nw}\mathbf{x}_3^{\rm nw})$, subject to $f_t(\dot{u}_1(r_1),\dot{u}_2(r_2),w;v)=0.$
          \item If $\rho_1$ or $\rho_2$ is abelian, then respectively
                $\mathbf{x}_2^{\rm nw}\in\{(\mathbf{x}_1^{\rm nw})^{\pm1}\}$ or $\mathbf{x}_2^{\rm nw}\in(\mathbf{x}_3^{\rm nw})^{\pm1}\}$.
          \item If $\rho_1,\rho_2$ are both NR, then $t$ is determined by $t=\kappa+\kappa^{-1}$ with
                $\Delta_{N(T_1)}(\kappa^2)=\Delta_{N(T_2)}(\kappa^2)=0$, and $\mathbf{x}_2^{\rm nw}=\mathbf{u}_{\kappa}(\eta)^\nu$, where $\nu\in\{\pm1\}$, and $\eta\ne 0,1$ can be arbitrary.
              %just to satisfy $\mathbf{x}_2^{\rm nw}\mathbf{x}_3^{\rm nw}\ne\mathbf{x}_3^{\rm nw}\mathbf{x}_2^{\rm nw}$.
        \end{itemize}
\end{enumerate}

In summary, the part of $\mathcal{X}^{\rm irr}(K)$ with $\mathbf{z}=\mathbf{e}$ has dimension $1$.

\subsection{$\mathbf{z}=\mathbf{p}$}\label{sec:g=p}

By Lemma \ref{lem:key} (b), $\mathbf{x}_i^{\rm nw}\in\mathsf{T}_+$, $i=1,2,3$, and $\mathbf{x}_i^{\rm ne}\in\mathsf{T}_+$, $i=1,2,4$. 
This forces $w\in\{2,t^2-2\}$. Up to finite ambiguity $r_3,r_4$ are determined by $t$ through $\dot{u}_3(r_3)=\dot{u}_4(r_4)=w$. 

Up to conjugacy we may assume $\mathbf{x}_1^{\rm nw}=\mathbf{d}(\kappa^{-1})$, $\mathbf{x}_1^{\rm ne}=\mathbf{u}_{\kappa}(1)$.
Then we can determine $\mathbf{x}_2^{\rm nw}$, $\mathbf{x}_2^{\rm ne}$ via $F_v^{T_1}$ as in (\ref{eq:Fv}), and further determine
$\mathbf{x}_3^{\rm nw}$, $\mathbf{x}_4^{\rm ne}$ via $F_v^{T_2}$.

Similarly as Case 2 in Section 4.2, $\rho$ is determined by $\mathbf{x}_3^{\rm ne}$ via ${\rm tr}(\mathbf{x}_3^{\rm ne}\mathbf{x}_3^{\rm nw})=u_3(r_3)$, ${\rm tr}(\mathbf{x}_3^{\rm ne}\mathbf{x}_4^{\rm ne})=t^2-u_4(r_4)$ and
$r_5:={\rm tr}(\mathbf{x}_3^{\rm nw}\mathbf{x}_4^{\rm ne}\mathbf{x}_3^{\rm ne})$, subject to
$$f_t\big(u_3(r_3),t^2-u_4(r_4),2;r_5\big)=0.$$
The irreducibility of $\rho$ is equivalent to $\mathbf{x}_3^{\rm ne}\notin\mathsf{T}_+$, which is equivalent to $\{u_3(r_3),u_4(r_4)\}\not\subseteq\{2,t^2-2\}$.

Thus, the part of $\mathcal{X}^{\rm irr}(K)$ contributed by $\mathbf{z}=\mathbf{p}$ has dimension $1$.

\subsection{$\mathbf{z}=\mathbf{d}(\kappa^2)$ with $\kappa+\kappa^{-1}=t\ne\pm2$}

Each of the cases $w\notin\{2,t^2-2\}$, $\mathbf{w}=\mathbf{e}$, $\mathbf{w}=\mathbf{p}$ is symmetric to one of the previous cases.
It remains to discuss the case $w=t^2-2$.

By Lemma \ref{lem:key} (d), for $i=1,2,3$, there exists $\xi_i$ such that
$\mathbf{x}_i^{\rm nw}$ equals $\mathbf{u}_\kappa(\xi_i)$ or $\mathbf{l}_\kappa(\xi_i)$.
Conjugating by 
$$\left(\begin{array}{cc} 0 & 1 \\ -1 & 0 \end{array}\right)$$ 
and/or replacing $\kappa$ by $\kappa^{-1}$ if necessary, up to conjugacy we may just assume $\mathbf{x}_1^{\rm nw}=\mathbf{u}_\kappa(1)$, so that $\mathbf{x}_1^{\rm ne}=\mathbf{u}_\kappa(-\kappa^{-2})$.

The possibility $\mathbf{x}_3^{\rm nw}=\mathbf{u}_\kappa(\xi)$ is ruled out by ${\rm tr}(\mathbf{x}_1^{\rm nw}\mathbf{x}_3^{\rm nw})=t^2-w=2$,
hence $\mathbf{x}_3^{\rm nw}=\mathbf{l}_{\kappa}(\zeta)$ for some $\zeta$, which is determined by
$$t^2-2=w={\rm tr}\big((\mathbf{x}_1^{\rm nw})^{-1}\mathbf{x}_3^{\rm nw}\big)
={\rm tr}\big(\mathbf{u}_{\kappa}(1)^{-1}\mathbf{l}_\kappa(\zeta)\big)=2-\zeta.$$

Up to finite ambiguity $r_1,r_2,r_3,r_4$ are determined by $t$, so are $\dot{u}_1(r_1)$, $\dot{u}_2(r_2)$, $u_3(r_3)$, $u_4(r_4)$.

If $\mathbf{x}_2^{\rm nw}=\mathbf{u}_{\kappa}(\eta)$, then $\dot{u}_1=2$, and $\eta$ is determined by
$$\dot{u}_2(r_2)={\rm tr}\big((\mathbf{x}_3^{\rm nw})^{-1}\mathbf{x}_2^{\rm nw}\big)
={\rm tr}\big(\mathbf{l}_{\kappa}(\zeta)^{-1}\mathbf{u}_\kappa(\eta)\big)=2-(4-t^2)\eta.$$
If $\mathbf{x}_2^{\rm nw}=\mathbf{l}_{\kappa}(\eta)$, then $\dot{u}_2=2$, and $\eta$ is determined by
$$\dot{u}_1(r_1)={\rm tr}\big((\mathbf{x}_1^{\rm nw})^{-1}\mathbf{x}_2^{\rm nw}\big)
={\rm tr}\big(\mathbf{u}_{\kappa}(1)^{-1}\mathbf{l}_\kappa(\eta)\big)=2-\eta.$$

Similarly, $\mathbf{x}_3^{\rm ne}$ has two possibilities.

In summary, the part of $\mathcal{X}^{\rm irr}(K)$ contributed by $\mathbf{z}=\mathbf{d}(\kappa^2)$ with $\kappa\ne\pm1$ has dimension $1$.

\subsection{Conclusion and discussion}\label{sec:conclusion}

Combining the results in Section 4.1--4.5 yields $\dim\mathcal{X}^{\rm irr}(K)=1$.
The proof of Theorem \ref{thm:dim=1} is finished by recalling that $\mathcal{X}^{\rm irr}(K)$ is Zariski open in $\mathcal{X}(K)$.

It is known that each ${\rm PSL}(2,\mathbb{C})$-representation of $K$ can be lifted to a ${\rm SL}(2,\mathbb{C})$-representation
(see \cite[Page 756]{BZ98}).
Thus, we have shown
\begin{cor}\label{cor:answer}
Under the assumptions {\rm(\ref{eq:assumption-0})}--{\rm(\ref{eq:assumption-2})}, the ${\rm PSL}(2,\mathbb{C})$-character variety of $K$ is $1$-dimensional.
\end{cor}

Recall the assumption $p_1,p_2,p_3,p_4\ge 3$. By \cite[Theorem 3.3]{Wu96}, $E_K$ contains a closed essential surface $\Sigma$ of genus $2$, and as claimed on \cite[Page 1]{Wu11}, $K$ is hyperbolic.
Thus, $K$ is a positive answer to Question \ref{que:BZ98}.

By \cite[Proposition 2.4]{CCGLS94}, if $K'$ is a knot with $\dim\mathcal{X}(K')>1$, then a closed incompressible surface in $E_{K'}$ is detected by ideal points by Culler-Shalen theory. Consequently, the surface $\Sigma\subset E_K$ is not detected by ideal points.

\begin{rmk}\label{rmk:related-work}
\rm Based on highly nontrivial computations, it was shown in \cite{CKT20} that in the exterior of each of the large hyperbolic knots $10_{152}, 10_{153}, 10_{154}$, no closed essential surface is detected by an ideal point of the character variety. As a consequence,
$\dim\mathcal{X}(L)=1$ for $L\in\{10_{152}, 10_{153}, 10_{154}\}$ (this can be verified by consulting \cite[Page 2]{CPY26}, where it is claimed that knots with $10$ crossings do not have high-dimensional character varieties, except for $10_{98},10_{99},10_{123}$).
Therefore, although the authors did not claim, they actually found answers to Question \ref{que:BZ98}.
\end{rmk}

By \cite[Theorem 3.6]{Wu96}, the Dehn filling $K(p/q)$ is Haken for any nontrivial slope $p/q\in\mathbb{Q}$.
This produces a lot of closed hyperbolic 3-manifolds containing an incompressible surface not detected by Culler-Shalen theory
(the first family of such manifolds was given in \cite[Section 10]{BZ98}), because by \cite[Theorem 5.8.2]{Th80}, $K(p/q)$ is hyperbolic except for finitely many $p/q$, and furthermore, $\dim\mathcal{X}(K(p/q))=0$ except for finitely many $p/q$.

To see the latter point, choose a meridian-longitude pair $(\mathfrak{m},\mathfrak{l})$ of $K$.
Let
$$\mathcal{R}^\ast(K)=\{\rho\in\mathcal{R}(K)\colon\rho(\mathfrak{m})=\mathbf{u}_{\kappa}(1)\ \text{for\ some\ }\kappa\ne 0\},$$
which is $1$-dimensional.
For each $\rho\in\mathcal{R}^\ast(K)$, note that $\rho(\mathfrak{l})\in\mathsf{T}_+$, as $\rho(\mathfrak{l})\rho(\mathfrak{m})=\rho(\mathfrak{m})\rho(\mathfrak{l})$; let $\kappa_\rho,\mu_\rho$ denote the upper-left entries of $\rho(\mathfrak{m}),\rho(\mathfrak{l})$, respectively.
Then $\rho\mapsto\kappa_\rho$ and $\rho\mapsto\mu_\rho$ are regular functions on $\mathcal{R}^\ast(K)$.
According to the knowledge on the A-polynomial (see \cite[Section 2]{CCGLS94}, and also \cite[Page 303]{LR03}), except for finitely many $p/q$, the condition $\kappa^p\mu^q=1$ cuts out a finite set.
So $\dim\mathcal{X}(K(p/q))=0$.

\bigskip

\noindent
Haimiao Chen (orcid: 0000-0001-8194-1264)\ \ \ \ {\it chenhm@math.pku.edu.cn} \\
Department of Mathematics, Beijing Technology and Business University, \\
Liangxiang Higher Education Park, Fangshan District, Beijing, China.

\end{document}